\documentclass[journal]{IEEEtran}
\usepackage{amsmath,amsfonts,amssymb,euscript, graphicx, 
epsfig,
enumerate,float,afterpage, subfigure, ifthen}%

\newtheorem{thm}{Theorem}

\newtheorem{lem}{Lemma}

\newcommand{\expect}[1]{\mathbb{E}\left\{#1\right\}}
\newcommand{\defequiv}{\mbox{\raisebox{-.3ex}{$\overset{\vartriangle}{=}$}}}

\newcommand{\norm}[1]{||{#1}||}

\newcommand{\bv}[1]{{\boldsymbol{#1} }}
\newcommand{\script}[1]{{{\cal{#1} }}}

\begin{document}

\title
  {Stochastic Optimization for Markov Modulated Networks 
  with Application to Delay Constrained Wireless Scheduling}
\author{Michael J. Neely , Sucha Supittayapornpong
\thanks{This material was presented in part at the $48$th 
IEEE Conf. on Decision and Control (CDC), 
Shanghai, China, Dec. 2009.} 
\thanks{The authors are with the  Electrical Engineering department at the University
of Southern California, Los Angles, CA.} 
\thanks{This material is supported in part  by one or more of 
the following: the DARPA IT-MANET program
grant W911NF-07-0028, the NSF Career grant CCF-0747525,  NSF grant 0540420, the 
Network Science Collaborative Technology Alliance sponsored
by the U.S. Army Research Laboratory W911NF-09-2-0053.}}

\markboth{}{Neely}

\maketitle

\begin{abstract} 
We consider a wireless system with a small
number of delay constrained users and a 
larger number of users without delay 
constraints.  We develop a scheduling
algorithm that reacts to time varying channels
and maximizes throughput utility (to within a desired
proximity), stabilizes
all queues, and satisfies the delay constraints.
The problem is solved by reducing the constrained
optimization to a set of weighted stochastic shortest path problems, 
which act as natural generalizations of max-weight policies
to Markov decision networks. 
We also present approximation results for the corresponding shortest
path problems, and  discuss the
additional complexity and delay incurred as compared to systems 
without delay constraints.  The solution technique is general and 
applies to other constrained 
stochastic decision problems. 
\end{abstract} 



\section{Introduction} 

This paper considers delay-aware scheduling in a multi-user wireless uplink or 
downlink with $K$ delay-constrained users and $N$ delay-unconstrained users, 
each with different transmission channels. 
The system operates in slotted time with 
normalized slots $t \in \{0, 1, 2, \ldots\}$.  Every slot, a random number of 
new packets arrive from each user.  
Packets are queued for eventual transmission, and every slot a scheduler 
looks at the queue backlog and the current channel states and chooses 
one channel to serve.  The number of packets  transmitted over that
channel depends on its current channel state. The goal is to stabilize all 
queues, satisfy average delay constraints for the delay-constrained
users, and drop  as few packets as possible. 

Without the delay constraints, this problem is a classical \emph{opportunistic scheduling} 
problem, and can be solved with  efficient max-weight algorithms based on Lyapunov 
drift and Lyapunov optimization (see \cite{now} and references therein).  The delay constraints
make the problem a much more complex \emph{Markov Decision Problem} (MDP).  While general 
methods for solving MDPs exist (see, for example, \cite{ross-prob}\cite{mdp-altman}\cite{meyn-control-book}), 
 they typically suffer from a curse of dimensionality.  Specifically, the number of queue state vectors grows
 exponentially in the number of queues.  Thus, a general problem with many queues has an intractably
 large state space. This creates non-polynomial implementation complexity for offline approaches
 such as linear programming \cite{ross-prob}\cite{mdp-altman}, and non-polynomial complexity and/or learning
 time for online or quasi online/offline approaches such as  $Q$-learning \cite{bertsekas-neural}\cite{q-learn-bert01}.

 We do not solve this fundamental curse of dimensionality. Rather, 
 we avoid this difficulty by focusing on the special structure that arises in a wireless 
 network with a \emph{relatively small number of delay-constrained users} (say, $K \leq 5$), 
 but with an arbitrarily large number of users without delay constraints (so that $N$ can be large). 
   This is an important scenario, particularly in cases when 
 the number of ``best effort'' users in a network is much larger than the number of 
 delay-constrained users.   We develop a solution that, on each slot, requires a computation that
 has a complexity that depends exponentially in $K$, but only polynomially
 in $N$.  Further, the resulting convergence times and delays are fully
 polynomial in the total number of queues $K+N$.   Our solution uses a concept of 
\emph{forced renewals} that introduces a deviation from optimality that can be made arbitrarily
small with a corresponding polynomial tradeoff in convergence time.  
Finally, we show that a simple Robbins-Monro iteration can be used to approximate the required
computations when channel and traffic statistics are unknown.    
Our methods are general and can be applied to other MDPs for 
networks with similar structure. 
 
 Related prior work on  delay optimality for multi-user opportunistic scheduling 
under special symmetric assumptions is developed in 
\cite{tass-server-allocation}\cite{edmund-thesis}\cite{anand-delay-it}, 
and single-queue delay optimization 
problems are treated in  
\cite{alvin-infocom2003}\cite{goyal-single-link}\cite{energy-delay-approx}\cite{q-learning-mimo} 
using dynamic programming and 
Markov Decision theory.  Approximate dynamic programming 
algorithms are applied to multi-queue switches in \cite{approx-dp-queue} and shown to perform
well in simulation. 
Optimal asymptotic energy-delay tradeoffs are developed for single queue systems in
\cite{berry-fading-delay}, and optimal energy-delay and utility-delay tradeoffs 
for multi-queue systems are treated in \cite{neely-energy-delay-it}\cite{neely-utility-delay-jsac}. 
The algorithms of \cite{neely-energy-delay-it}\cite{neely-utility-delay-jsac} 
have very low complexity and provably
converge quickly even for large networks, although the tradeoff-optimal delay guarantees 
they achieve do not necessarily optimize the coefficient multiplier in the delay expression. 

Our approach in the present paper treats the MDP problem associated with
delay constraints using 
Lyapunov drift and Lyapunov optimization theory \cite{now}.  This theory has been used to 
stabilize queueing networks \cite{tass-server-allocation} and provide utility optimization \cite{neely-thesis}\cite{stolyar-greedy}\cite{atilla-fairness-ton}\cite{atilla-jsac}\cite{now} via simple \emph{max-weight} principles.
We extend the max-weight
principles to treat networks with  \emph{Markov decisions}, 
where the network costs depend on both the control actions taken and the current
state (such as the queue state) the system is in.  For each cost constraint  
we define a \emph{virtual queue}, and show that the constrained
MDP can be solved using Lyapunov drift theory implemented over a variable-length
frame, where ``max-weight'' rules are replaced with weighted stochastic shortest path problems. 
This is similar to the Lagrange multiplier approaches used in the related works 
\cite{energy-delay-approx}\cite{q-learning-mimo} that treat 
power minimization for single-queue wireless links with an average delay constraint. 
The work in  \cite{energy-delay-approx}
uses stochastic approximation with a 2-timescale argument and a limiting
ordinary differential equation.   The work in \cite{q-learning-mimo} 
treats a single-queue MIMO system using primal-dual updates \cite{bertsekas-nonlinear}. 
Our virtual queues are similar to the Lagrange Multiplier updates in \cite{energy-delay-approx}\cite{q-learning-mimo}.  
However, we treat multi-queue systems, and we use a different 
analytical approach that emphasizes
stochastic shortest paths over variable length frames.  Because of this, 
our approach can be 
used in conjunction with a variety of existing techniques for solving shortest 
path problems (see, for example, 
\cite{bertsekas-neural}). 
We use a Robbins-Monro technique that is adapted to this context, together
with a \emph{delayed queue analysis} to uncorrelate past samples from current queue 
states.
Our resulting algorithm has an implementation complexity that 
grows exponentially in the number of delay-constrained queues $K$, but polynomially
in the number of delay-unconstrained queues $N$.  Further, we obtain
polynomial bounds  on convergence times and delays.

The next section describes the network model. 
Section \ref{section:alg} 
presents  the weighted stochastic shortest-path algorithm.
 Section \ref{section:approx} describes approximate implementations, and Section \ref{section:simulation} 
 presents a simple simulation example. 

 \section{Network Model} \label{section:model} 
 
 Consider a wireless queueing network that operates in discrete time with timeslots $t \in \{0, 1, 2,\ldots\}$. 
 The network has $K$ \emph{delay-constrained queues} and $N$ \emph{stability-constrained queues}, 
 for a total of $K+N$ queues indexed by sets $\script{K} \defequiv \{1, \ldots, K\}$ and  
 $\script{N} \defequiv \{K+1, \ldots, K+N\}$. 
 The queues store fixed-length packets for
 transmission over their wireless channels.  
 Every timeslot, new packets randomly arrive to each queue, and we let $\bv{A}(t) = (A_1(t), \ldots, A_{K+N}(t))$
 represent the random packet arrival vector.   The stability-constrained queues
 have an infinite buffer space.  The delay-constrained queues have a finite buffer space that can store $b$ packets
 (for some positive integer $b$).  
 The network channels can vary from slot to slot, and we let $\bv{S}(t) = (S_1(t), \ldots, S_{K+N}(t))$
 be the channel state vector on slot $t$, representing conditions that affect transmission rates. We assume
 the stacked vector  $[\bv{A}(t), \bv{S}(t)]$ is independent and identically distributed (i.i.d.) over slots, with possibly correlated entries
 on the same slot. 
 
  Every slot $t$, the network controller observes the channel states $\bv{S}(t)$ and chooses 
  a \emph{transmission rate vector} $\bv{\mu}(t) = (\mu_1(t), \ldots, \mu_{K+N}(t))$, being a vector of non-negative integers.
  The choice of $\bv{\mu}(t)$ is constrained to a set $\Gamma_{\bv{S}(t)}$ that depends on the current $\bv{S}(t)$.  A simple
  example is a system with ON/OFF channels where the controller can transmit a single packet
  over at most one ON channel per slot, as in \cite{tass-server-allocation}.  In this example, $\bv{S}(t)$ is a binary vector
  of channel states, and $\Gamma_{\bv{S}(t)}$ restricts $\bv{\mu}(t)$ to be a binary vector with at most one non-zero entry and
  with  $\mu_i(t) = 0$ whenever $S_i(t) = 0$.  We assume that for each possible channel state vector 
  $\bv{S}$, the set $\Gamma_{\bv{S}}$ has the property
  that for any $\bv{\mu} \in \Gamma_{\bv{S}}$, the vector $\bv{\mu}'$ is also in $\Gamma_{\bv{S}}$, where $\bv{\mu}'$ is formed
  from $\bv{\mu}$ by setting one or more entries to $0$. 
  In addition to constraining $\bv{\mu}(t)$ to take values in $\Gamma_{\bv{S}(t)}$ every slot $t$, we shall soon 
  also restrict the $\mu_k(t)$ values
  for the delay-constrained  queues $k \in \script{K}$ to be at most the current number of packets in queue $k$.  This is a natural restriction,
  although we \emph{do not} place such a restriction on the stability-constrained queues $n \in \script{N}$.  This is a technical detail that will be
  important later, when we show that the \emph{effective dimension} of the resulting Markov decision problem is $K$, independent
  of the number of stability-constrained queues $N$. 
  
Let $\bv{Q}(t) = (Q_1(t), \ldots, Q_{K+N}(t))$ represent the vector of current queue backlogs, and define 
$d_n(t) = A_n(t) - \mu_n(t)$.  The queue dynamics for the 
stability-constrained queues are:\footnote{For simplicity of exposition later, we have allowed the stability-constrained 
queues $Q_n(t)$  to serve
newly arriving data.  This can be modified easily by introducing a delay by one slot, so that the ``new arrivals'' to the stability-constrained
queues actually arrived one slot ago.}
\begin{equation} \label{eq:q-update} 
Q_n(t+1) = \max[Q_n(t) + d_n(t), 0] \: \: \forall n \in \script{N} 
\end{equation} 
where the $\max[\cdot,0]$ operation allows, in principle, a service variable $\mu_n(t)$ to be independent of whether or not
$Q_n(t)$ is empty.  

The delay-constrained queues have a different queue dynamic.  Because of the finite buffer, we 
must allow packet dropping.  Let $D_k(t)$ be the number of dropped packets on slot $t$.  The queue dynamics for the delay-constrained
queues are given by:  
\begin{equation} \label{eq:z-update} 
 Q_k(t+1) = Q_k(t) - \mu_k(t) - D_k(t) + A_k(t) \: \: \forall k \in \script{K} 
 \end{equation} 
Note that this does not have any $\max[\cdot, 0]$ operation, because we will force the $\mu_k(t)$ and $D_k(t)$ decisions to be 
such that we never serve or drop packets that we do not have.  The precise constraints on these decision variables is given 
after the introduction of a \emph{forced renewal event}, defined in the next subsection. 

\subsection{Forced Renewals} 

To force the delay-constrained queues to repeatedly visit a \emph{renewal state} of being simultaneously empty, 
at the end of every slot, with probability $\phi>0$ we independently drop all unserved packets in 
all delay constrained queues $k \in \script{K}$.  
The stability-constrained queues do not experience such forced drops. 
Specifically, let 
$\phi(t)$ be an i.i.d. Bernoulli process that is $1$ with probability $\phi$ every slot $t$, and $0$ otherwise.  Assume 
$\phi(t)$ is independent of $[\bv{A}(t), \bv{S}(t)]$.  If $\phi(t) =1$, we say slot $t$ experiences a \emph{forced renewal event}. 
The decision options for $\mu_k(t)$ and $D_k(t)$ for $k \in \script{K}$ are then additionally constrained as follows:  If $\phi(t) = 0$, then:
\begin{eqnarray*}
\mu_k(t) &\in& \{0, 1, \ldots, Q_k(t)\} \\
D_k(t) &\in& \{\max[A_k(t)+Q_k(t) - b, 0], \ldots, A_k(t)\} 
\end{eqnarray*}
so that during normal operation,  we can serve at most $Q_k(t)$  packets from queue $k$ (so new arrivals cannot be served),
and we can drop only new arrivals, necessarily dropping any new arrivals that would exceed the finite buffer capacity. 
However, if $\phi(t) = 1$ we have: 
\begin{eqnarray*}
\mu_k(t) &\in& \{0, 1, \ldots, Q_k(t)\} \\
D_k(t) &=& Q_k(t) - \mu_k(t) + A_k(t) 
\end{eqnarray*}
So that $\mu_k(t)$ is constrained as before, but $D_k(t)$ is then equal to the remaining packets (if any) at the end of the slot.

We shall optimize the system under the assumption that the forced renewal process $\phi(t)$ is uncontrollable.  This provides
an analyzable system that lends itself to simple approximations, as shown in later parts of the paper. While these forced renewals
create inefficiency in the system, the rate of dropped packets due to forced renewals is at most $(Kb + \sum_{k=1}^K\expect{A_k(t)})\phi$, 
which assumes the worst case of dropping a full buffer plus all new arrivals every renewal event. 
This value can be made arbitrarily small with a small choice of $\phi$. For problems such as minimizing the average drop rate subject to 
delay constraints in the delay-constrained queues and stability in the stability-constrained queues, it can be shown that 
this $O(\phi)$ term bounds the gap between system optimality without forced renewals and system optimality with forced renewals.  
Formally, this can be shown by a simple sample path argument:  A system optimized without forced renewals has a performance that is no better
than a system with forced renewals, but where all ``drops'' from forced renewals are counted as delivered throughput, and where all other decisions mimic those of the prior system.  We omit a formal argument for brevity.   In Theorem \ref{thm:1} we show the disadvantage of using 
a small value of $\phi$ is that our average queue bounds for the stability-constrained queues is $O(1/\phi)$. 

Define a \emph{renewal frame} as the sequence of slots starting just after a renewal event and ending at the next renewal event. 
Assume that $\bv{Z}(0) =\bv{0}$, so that time $0$ starts the first renewal frame.  Define $t_0=0$, and let $t_r$ and for $r \in \{1, 2, \ldots\}$
represent the sequence that marks the beginning of each renewal frame.  
For $r \in \{0, 1, 2, \ldots\}$, define $T_r \defequiv t_{r+1}-t_r$ as the duration
of the $r$th renewal frame. Note that $\{T_r\}_{r=0}^{\infty}$ are  i.i.d. geometric random variables with $\expect{T_r} = 1/\phi$.

\subsection{Markov Decision Notation} \label{section:notation} 

Define $\omega(t) \defequiv [\bv{A}(t), \bv{S}(t)]$ as the observed arrivals and channels of the network on slot $t$, 
and define the random network event $\Omega(t) \defequiv [\omega(t), \phi(t)]$.  Then $\Omega(t)$ is i.i.d.
over slots.  We can summarize the control decision constraints of the previous section with the following simple notation: 
Let $\script{Z} \defequiv \{0, 1, \ldots, b\}^K$ be the $K$-dimensional state space for the delay-constrained queues, and 
let $z(t) \defequiv (Q_k(t))_{k\in\script{K}}$ represent the current state of these queues.  Every slot $t$, the controller observes
the random event 
$\Omega(t)$ and the queue state $z(t)$, and makes a \emph{control action} $\alpha(t)$, which determines all decision variables
$\mu_i(t)$ for $i \in \{1, \ldots, K+N\}$ and $D_k(t)$ for $k \in \script{K}$,  chosen in a set $\script{A}_{\Omega(t), z(t)}$ 
that depends on $\Omega(t)$ and $z(t)$.  Note that, indeed, all of our decision variables as described in the previous subsection 
are constrained only
in terms of $\Omega(t)$ and $z(t)$, and in particular the queue states $Q_n(t)$ for $n \in \script{N}$ do not constrain our decisions. 

Recall that $d_n(t) \defequiv A_n(t) - \mu_n(t)$. The $\alpha(t)$, $\Omega(t)$, $z(t)$ together affect the vector $\bv{d}(t) = (d_n(t))_{n\in\script{N}}$
through a 
deterministic function $\hat{d}_n(\alpha(t), \Omega(t), z(t))$: 
\begin{equation} \label{eq:dn} 
 d_n(t) = \hat{d}_n(\alpha(t), \Omega(t), z(t)) \: \: \forall n \in \script{N} 
 \end{equation} 
Further, $\alpha(t)$, $\Omega(t)$, $z(t)$ together define the \emph{transition probabilities} 
from $z(t)$ to $z(t+1)$, defined for all states $i$ and $j$ in $\script{Z}$: 
\begin{equation} \label{eq:pij} 
 P_{ij}(\alpha, \Omega) = Pr[z(t+1) = j|z(t) = i, \alpha(t) = \alpha, \Omega(t) = \Omega] 
 \end{equation} 
From the equation \eqref{eq:z-update} we find that $P_{ij}(\alpha, \Omega) \in \{0,1\}$, 
so that next states $z(t+1)$ are deterministic
given $\alpha(t)$, $\Omega(t)$, $z(t)$. Finally, we define a general \emph{penalty vector} $\bv{y}(t) = (y_0(t), y_1(t), \ldots, y_L(t))$, 
for some integer $L\geq 0$, 
where penalties $y_l(t)$ are deterministic functions of $\alpha(t)$, $\Omega(t)$, $z(t)$: 
\begin{equation} \label{eq:yl} 
y_l(t) \defequiv \hat{y}_l(\alpha(t), \Omega(t), z(t))
\end{equation} 
 For example, penalty $y_0(t)$ can be defined as the total number of dropped packets on slot $t$ by 
 defining $y_0(t) = \sum_{k\in\script{K}} D_k(t)$, which is indeed a function of $\alpha(t)$, $\Omega(t)$, $z(t)$. 

We assume throughout that all of the above deterministic functions are bounded, so that there is a finite constant $\beta$ such that
for all $l \in \{0, 1, \ldots, L\}$, all $n \in \script{N}$, and all slots $t$ we have: 
\begin{equation} \label{eq:boundedness} 
 |y_l(t)| \leq \beta \: \: , \: \: |d_n(t)| \leq \beta 
 \end{equation}

\subsection{The Optimization Problems} 

A control policy is a method for choosing actions $\alpha(t) \in \script{A}_{\Omega(t), z(t)}$ over slots $t \in \{0, 1, 2, \ldots\}$. 
We restrict to causal policies that make decisions with knowledge of the past but without knowledge of the future.  Suppose
a particular control policy is given.  Define time averages $\overline{Q}_n$ and $\overline{y}_l$ for $n\in\script{N}$ and $l \in \{0, 1, \ldots, L\}$
by: 
\begin{eqnarray*}
\overline{Q}_n &\defequiv& \limsup_{t\rightarrow\infty}\frac{1}{t}\sum_{\tau=0}^{t-1} \expect{Q_n(\tau)} \\
\overline{y}_l &\defequiv& \limsup_{t\rightarrow\infty} \frac{1}{t}\sum_{\tau=0}^{t-1} \expect{y_l(\tau)} 
\end{eqnarray*}
Our goal is to design a control policy to solve the following stochastic optimization problem: 
\begin{eqnarray} 
\mbox{Minimize:} & \overline{y}_0 \label{eq:p0} \\
\mbox{Subject to:} & \overline{y}_l \leq 0 \: \: \forall l \in \{1, \ldots, L\} \label{eq:p1}  \\
& \overline{Q}_n < \infty \: \: \forall n \in \script{N} \label{eq:p2} \\
& \alpha(t) \in \script{A}_{\Omega(t), z(t)} \: \: \forall t \in \{0, 1, 2, \ldots\} \label{eq:p3} 
\end{eqnarray} 

That is, we desire to minimize the time average of the $y_0(t)$ penalty, subject to time average constraints on the 
other penalties, and subject to queue stability (called \emph{strong stability}) 
for all stability-constrained queues. The general structure \eqref{eq:p0}-\eqref{eq:p3} fits a variety of
network optimization problems.  For example, if we define $y_0(t)$ as the sum packet drops $\sum_{k\in\script{K}} D_k(t)$, 
define $L=K$, and define $y_k(t) = Q_k(t) - Q_{av}$ for all $k \in \script{K}$ (for some positive constant $Q_{av}$), 
then the problem \eqref{eq:p0}-\eqref{eq:p3} seeks to 
minimize the total packet drop rate, subject to an average backlog of at most $Q_{av}$ in all delay-constrained queues 
$k \in \script{K}$, and subject to stability of all stability-constrained queues $n \in \script{N}$. 

Alternatively, to enforce an average \emph{delay} constraint $W_{av}$ at all queues $k \in \script{K}$
(for some positive number $W_{av}$), we can 
define penalties: 
\[ y_k(t) = Q_k(t) - (A_k(t) - D_k(t))W_{av} \: \: \: \forall k \in \script{K}  \]
Note that the time average of $(A_k(t) - D_k(t))$ is the number 
$\tilde{\lambda}_k$, the average arrival rate of (non-dropped) packets to queue $k$.  Hence, the constraint $\overline{y}_k \leq 0$
is equivalent to: 
\[ \overline{Q}_k - \tilde{\lambda}_kW_{av} \leq 0 \]
However, by Little's theorem \cite{bertsekas-data-nets} we have $\overline{Q}_k = \tilde{\lambda}_k\overline{W}_k$, where $\overline{W}_k$
is the average delay for queue $k$, and so the 
constraint $\overline{y}_k \leq 0$ ensures $\overline{W}_k \leq W_{av}$ (assuming $\lambda_k >0$).

In the following, we develop a dynamic algorithm that can come arbitrarily close to solving
the problem \eqref{eq:p0}-\eqref{eq:p3}.  Our solution is general and applies to any other discrete time 
Markov decision problem
on a general finite state space $\script{Z}$, random events $\Omega(t) = [\omega(t), \phi(t)]$ (for forced
renewal process $\phi(t)$),
control actions $\alpha(t)$ in a general set $\script{A}_{\Omega(t), z(t)}$, queue equations \eqref{eq:q-update} with 
$d_n(t)$ given in the form \eqref{eq:dn}, transition probabilities in the form \eqref{eq:pij}, and penalties in the form
\eqref{eq:yl}.

\subsection{Slackness Assumptions}  \label{section:slackness} 

Suppose the problem \eqref{eq:p0}-\eqref{eq:p3} is \emph{feasible}, so that there exists a policy that satisfies the constraints.
It can be shown that the 
constraint $\overline{Q}_n < \infty$ implies that $\overline{d}_n \leq 0$  \cite{sno-text}, and so the following modified problem is feasible
whenever the original one is: 
\begin{eqnarray}
\mbox{Minimize:} & \overline{y}_0 \label{eq:p0m} \\
\mbox{Subject to:} & \overline{y}_l \leq 0 \: \: \forall l \in \{1, \ldots, L\} \label{eq:p1m}  \\
& \overline{d}_n \leq 0  \: \: \forall n \in \script{N} \label{eq:p2m} \\
& \alpha(t) \in \script{A}_{\Omega(t), z(t)} \: \: \forall t \in \{0, 1, 2, \ldots\} \label{eq:p3m} 
\end{eqnarray}
Define $y_0^{opt}$ as the infimum of $\overline{y}_0$ for the problem \eqref{eq:p0m}-\eqref{eq:p3m}, necessarily
being less than or equal to the corresponding infimum of the original problem \eqref{eq:p0}-\eqref{eq:p3}.\footnote{Recall
that $y_0^{opt}$ is defined assuming forced renewals of probability $\phi$.  Thus, $y_0^{opt}$ is within a gap of $O(\phi)$ of the 
minimum cost without such forced renewals.}  We show
in Theorem \ref{thm:1} that, under a suitable slackness condition, the value of $y_0^{opt}$ can be approached arbitrarily
closely while maintaining $\overline{Q}_n < \infty$ for all queues $n \in \script{N}$.  Thus, under that slackness condition, 
 $y_0^{opt}$ is also the infimum of $\overline{y}_0$
for the original problem \eqref{eq:p0}-\eqref{eq:p3}.

The problem \eqref{eq:p0m}-\eqref{eq:p3m} is a constrained Markov decision problem (MDP) with state $(\Omega(t), z(t))$. Under
mild assumptions (such as this state space being finite, and the action space $\script{A}_{\Omega, z}$ being finite for each 
$(\Omega, z)$) the MDP has an \emph{optimal stationary policy} that 
chooses actions $\alpha(t) \in \script{A}_{\Omega(t), z(t)}$ every slot $t$
as a stationary and possibly randomized function of the state
$(\Omega(t), z(t))$ only.  We call such policies \emph{$(\Omega,z)$-only policies}.  Because this system experiences regular renewals,
the performance of any $(\Omega, z)$-only policy can be characterized by ratios of expectations over one renewal frame. 
Thus, we make the following assumption. 

\emph{Assumption 1:}  There is an $(\Omega, z)$-only policy $\alpha_1^*(t)$ that satisfies the following over any renewal frame: 
\begin{eqnarray}
\frac{\expect{\sum_{\tau=t_r}^{t_r+T_r-1} y_0^*(\tau)}}{1/\phi} &=& y_0^{opt} \label{eq:a1-1} \\
\frac{\expect{\sum_{\tau=t_r}^{t_r + T_r-1} d_n^*(\tau)}}{1/\phi} &\leq& 0 \: \: \forall n \in \script{N} \label{eq:a1-2} \\
\frac{\expect{\sum_{\tau=t_r}^{t_r+T_r-1}y_l^*(\tau)}}{1/\phi} &\leq& 0 \: \: \forall l \in \{1, \ldots, L\} \label{eq:a1-3} 
\end{eqnarray} 
where $T_r$ is the size of the renewal frame, with $\expect{T_r} = 1/\phi$, and $y_l^*(\tau)$, $d_n^*(\tau)$ are values under
the policy $\alpha^*(t)$ on slot $\tau$ of the renewal frame. 

We emphasize that Assumption 1 is mild and holds whenever the problem \eqref{eq:p0m}-\eqref{eq:p3m} is feasible and has an
optimal stationary policy (i.e., an optimal $(\Omega, z)$-only policy).  
We now make the following stronger assumption that there exists an $(\Omega,z)$-only policy that can meet the constraints 
\eqref{eq:a1-2}-\eqref{eq:a1-3} with ``$\epsilon$-slackness,'' without caring what average value of $y_0(t)$ this policy generates.
This assumption is related to standard ``Slater-type'' assumptions in optimization theory \cite{bertsekas-nonlinear}. 

\emph{Assumption 2:} There is a value $\epsilon>0$ and an $(\Omega, z)$-only policy $\alpha_2^*(t)$ (typically different from 
policy $\alpha_1^*(t)$ in Assumption 1) 
that satisfies the following over any renewal frame: 
\begin{eqnarray}
\frac{\expect{\sum_{\tau=t_r}^{t_r + T_r-1} d_n^*(\tau)}}{1/\phi} &\leq& -\epsilon \: \: \forall n \in \script{N} \label{eq:a2-1} \\
\frac{\expect{\sum_{\tau=t_r}^{t_r+T_r-1}y_l^*(\tau)}}{1/\phi} &\leq& -\epsilon \: \: \forall l \in \{1, \ldots, L\} \label{eq:a2-2} 
\end{eqnarray} 

We show in Theorem \ref{thm:1} that systems that satisfy Assumption 2 with larger values of $\epsilon$ can operate with
smaller average queue sizes in the stability-constrained queues.

 \section{The Dynamic Control Algorithm} \label{section:alg} 
 
 To solve the problem \eqref{eq:p0}-\eqref{eq:p3}, we
 extend the framework of \cite{now} to a case of variable length frames.  Specifically, for each 
 of the $L$ penalty constraints $\overline{y}_l \leq 0$, we define a
 \emph{virtual queue} $X_m(t)$ that is initialized to zero and that has dynamic update equation: 
 \begin{equation} \label{eq:x-update} 
 X_l(t+1) = \max[X_l(t) + y_l(t), 0] \: \: \forall l \in \{1, \ldots, L\}
 \end{equation} 
 where $y_l(t) = \hat{y}_{l}(\alpha(t), \Omega(t), z(t))$ is the $l$th penalty incurred on slot $t$ by  a particular 
 action $\alpha(t) \in \script{A}_{\Omega(t), z(t)}$.    The intuition is that if 
 the virtual queue $X_l(t)$ is stable, then the time average of $y_l(t)$ must be non-positive.  This turns the time 
 average constraint into a simple queue stability problem.
  
 \subsection{Lyapunov Drift} 
 
 Define $\bv{X}(t)$ as a vector of all virtual queues $X_l(t)$ for $l \in \{1, \ldots, L\}$.
 Define  $\bv{\Theta}(t)$ as the combined vector of all virtual queues and all stability-constrained queues: 
 \[ \bv{\Theta}(t) \defequiv [\bv{X}(t), (Q_n(t))_{n\in\script{N}}] \]
 Assume all queues are initially empty, so that $\bv{\Theta}(0) = \bv{0}$. 
 Define the following quadratic function: 
 \[ L(t) \defequiv \frac{1}{2}\sum_{n \in \script{N}} Q_n(t)^2 + \frac{1}{2}\sum_{l=1}^L X_l(t)^2 \]
Let $t_r$ be the start of a renewal frame, with duration $T_r$. Define the 
\emph{frame-based conditional Lyapunov drift} 
$\Delta(t_r)$ as follows:
\begin{eqnarray} 
 \Delta(t_r) \defequiv  \expect{L(t_r+T_r) - L(t_r)\left|\right. \bv{\Theta}(t_r), z(t_r) = 0}\label{eq:drift-def} 
 \end{eqnarray} 
 Note that $\Delta(t_r)$ is a function of the initial state $\bv{\Theta}(t_r)$ and the policy implemented during the 
 frame, where expectations are with respect to the random 
 events that can take place and the possibly random control actions made. 
 The explicit conditioning on $z(t_r)=0$ in \eqref{eq:drift-def} will be suppressed in the remainder
 of this paper, as this conditioning is implied given that $t_r$ starts a renewal frame. 
 
  It is important to note the 
 following subtlety:  The implemented policy $\alpha(t)$ may not be stationary and/or may depend on the queue
values $\bv{Q}(t)$ (which can be different on each renewal interval), 
and so actual system events are 
not necessarily i.i.d. over different
 renewal frames.   However, these frames are useful because we will analytically compare
 the Lyapunov drift of the actual implemented policy over a frame  to the corresponding drifts of the $(\Omega,z)$-only 
 policies of Assumptions 1 and 2.

 \begin{lem} \label{lem:drift} (Lyapunov Drift) Under any network control policy that chooses $\alpha(\tau) \in \script{A}_{\Omega(\tau), z(\tau)}$
 for all slots $\tau$ during a renewal frame $\tau \in \{t_r, \ldots, t_r +T_r-1\}$,  and for any initial queue values 
 $\bv{\Theta}(t_r)$, we have:
 \begin{eqnarray}
 \Delta(t_r) &\leq& B/\phi^2  + \expect{D(\bv{\Theta}(t_r))|\bv{\Theta}(t_r)} \label{eq:rhs-drift} 
 \end{eqnarray} 
 where $D(\bv{\Theta}(t_r))$ is defined: 
 \begin{eqnarray} 
 D(\bv{\Theta}(t_r)) &\defequiv& 
   \sum_{n \in \script{N}} Q_n(t_r)\sum_{\tau=t_r}^{t_r+T_r-1} d_n(\tau)\nonumber  \\
  && + \sum_{l=1}^L X_l(t_r)\sum_{\tau=t_r}^{t_r+T_r-1} y_l(\tau) \label{eq:D-def}
 \end{eqnarray}
 and where $B$ is a finite constant defined: 
 \begin{eqnarray*}
 B \defequiv \frac{(2-\phi)\beta^2(N+L)}{2}
 \end{eqnarray*}
 where we recall $\beta$ is the bound in \eqref{eq:boundedness}.
 \end{lem}

 \begin{proof} 
 For any $l \in \{1, \ldots, L\}$ and any $\tau \in \{t_r, \ldots, t_r+T_r-1\}$ we have by squaring \eqref{eq:x-update}:
 \begin{eqnarray*}
  X_l(\tau+1)^2 &\leq& (X_l(\tau) + y_l(\tau))^2 \\
  &=& X_l(\tau)^2 + y_l(\tau)^2 + 2X_l(\tau)y_l(\tau) \\
  &=& X_l(\tau)^2 + y_l(\tau)^2 + 2X_l(t_r)y_l(\tau)  \\
  && + 2[X_l(\tau) -X_l(t_r)]y_l(\tau) \\
  &\leq& X_l(\tau)^2 + \beta^2 + 2X_l(t_r)y_l(\tau) + 2\beta^2(\tau-t_r)
  \end{eqnarray*}
  where the final inequality holds because the change in $X_l(\tau)$ on any slot is at most $\beta$, 
  as is the magnitude of $y_l(\tau)$.  Summing the above over $\tau \in \{t_r, \ldots, t_r+T_r-1\}$
  and dividing by $2$ yields: 
  \begin{eqnarray}
  \frac{X_l(t_r+T_r)^2-X_l(t_r)^2}{2} &\leq& \frac{T_r\beta^2 + \beta^2T_r(T_r-1)}{2} \nonumber \\
  && + X_l(t_r)\sum_{\tau=t_r}^{t_r+T_r-1}y_l(\tau) \label{eq:123} \\
  &=& \frac{\beta^2T_r^2}{2} + X_l(t_r)\sum_{\tau=t_r}^{t_r+T_r-1}y_l(\tau) \nonumber \\
  && \label{eq:driftlem1} 
  \end{eqnarray}
  where \eqref{eq:123}  uses the identity: 
  \[ \sum_{\tau=t_r}^{t_r+T_r-1}(\tau-t_r) = T_r(T_r-1)/2 \]
  Similarly, it can be shown for any $n \in \script{N}$: 
  \begin{eqnarray}
  \frac{Q_n(t_r+T_r)^2 - Q_n(t_r)^2}{2} &\leq&  \frac{\beta^2T_r^2}{2}  \nonumber \\
  && \hspace{-.3in}+ Q_n(t_r)\sum_{\tau=t_r}^{t_r+T_r-1}d_n(\tau) \label{eq:driftlem2}
  \end{eqnarray}
  Summing \eqref{eq:driftlem1} and \eqref{eq:driftlem2} over $l \in \{1, \ldots, L\}$, $n\in\script{N}$, 
  taking conditional expectations, and noting that the second moment of a geometric random variable $T_r$ with success probability 
  $\phi$ is given by $(2-\phi)/\phi^2$ proves the result.  
 \end{proof} 
 
 \subsection{The Frame-Based Drift-Plus-Penalty Algorithm} \label{section:dpp-alg} 

Let $V\geq0$ be a 
 non-negative parameter 
 that we use to affect proximity to the optimal solution.
 Our dynamic algorithm initializes all virtual and actual queue states to 0, and designates $t_0=0$ as the start
 of the first renewal frame.  Then: 
 \begin{itemize} 
 \item For each frame $r \in \{0, 1, 2, \ldots\}$, observe the vector of virtual and actual queues 
 $\bv{\Theta}(t_r)$ and implement a policy over the course of the frame to minimize the following ``drift-plus-penalty''
 expression: 
 \begin{equation} \label{eq:dpp} 
 \expect{ D(\bv{\Theta}(t_r)) + V\sum_{\tau=t_r}^{t_r+T_r-1} y_0(\tau) \left|\right.\bv{\Theta}(t_r)} 
 \end{equation} 
 \item During the course of the frame, update virtual and actual queues every slot
 by \eqref{eq:q-update} and \eqref{eq:x-update}, and update state $z(t)$ by \eqref{eq:pij}.  
 At the end of the frame, go back to the preceding step. 
 \end{itemize} 

The decision rule \eqref{eq:dpp} generalizes the drift-plus-penalty rule in \cite{now}\cite{neely-energy-it} to 
a variable frame system.  The problem of designing a policy to minimize \eqref{eq:dpp} 
is a   \emph{weighted
stochastic shortest path problem}, where weights are virtual and actual queue backlogs at the start of the frame. 
Finding such a policy is non-trivial, and often can only be done in an approximate context.  In the next sub-section,
we present the performance of the algorithm, under the assumption that we have an algorithm to approximate
\eqref{eq:dpp}.   In Section \ref{section:approx} we consider various such approximation methods.

  \subsection{Performance Theorem} 
 
 For constants $C\geq 0$, $\delta\geq0$, define a \emph{$(C,\delta)$-approximation} of \eqref{eq:dpp} to be a policy for choosing 
 $\alpha(t)$ over a frame (consisting of slots $\tau\in\{t_r, \ldots, t_r+T_r-1\}$)  that yields 
 a total drift-plus-penalty that is less than or equal to that of any other policy, plus an error term parameterized by $C$ and $\delta$:
 \begin{eqnarray} 
 \expect{D(\bv{\Theta}(t_r)) + V\sum_{\tau=t_r}^{t_r+T_r-1} y_0(\tau)\left|\right.\bv{\Theta}(t_r)} 
\leq  \nonumber \\
\expect{D^*(\bv{\Theta}(t_r)) +  V\sum_{\tau=t_r}^{t_r+T_r-1} y_0^*(\tau)\left|\right.\bv{\Theta}(t_r)} \nonumber \\
+  C + \delta\sum_{n\in\script{N}}Q_n(t_r) + \delta\sum_{l=1}^LX_l(t_r) +  V\delta \label{eq:delta-opt}
 \end{eqnarray}
 where $D^*(\bv{\Theta}(t_r))$ and $y^*_0(\tau)$ represent \eqref{eq:D-def} and \eqref{eq:yl}, respectively, under
 any alternative algorithm $\alpha^*(t)$ that can be implemented during the slots $\tau \in \{t_r, \ldots, t_r+T_r-1\}$ of the frame. 
  Note that an exact minimization of the stochastic shortest path problem \eqref{eq:dpp} is a $(C, \delta)$-approximation for $C=\delta=0$. 
  
  \begin{thm} \label{thm:1} Suppose 
  Assumptions 1 and 2 hold for
  a given $\epsilon>0$.  Fix $V\geq 0$, $C\geq0$, $\delta\geq0$, 
  and suppose we use a $(C,\delta)$-approximation every frame.  If $\epsilon > \phi\delta$, then all virtual and actual 
  queues are strongly stable, and so all desired constraints \eqref{eq:p1}-\eqref{eq:p3} are satisfied.   In particular, for all positive 
  integers $R$, the average queue sizes satisfy: 
  \begin{eqnarray}
  \frac{1}{R}\sum_{r=0}^{R-1} \left[\sum_{n \in \script{N}}\expect{Q_n(t_r)} + \sum_{l=1}^L \expect{X_l(t_r)}\right]\leq \nonumber \\
  \frac{B/\phi +  C\phi  + V(\phi\delta + 2\beta)}{\epsilon - \phi\delta}   \label{eq:opt-thm1} 
  \end{eqnarray}
  Further, the time average penalty satisfies: 
  \begin{eqnarray}
 &\hspace{-.5in}  \limsup_{t\rightarrow\infty} \frac{1}{t}\sum_{\tau=0}^{t-1} \expect{y_0(\tau)} \leq  y_0^{opt} + \nonumber \\
& \hspace{-.2in}   \frac{B/\phi + C\phi}{V}   + \phi\delta[1 + (\beta - y_0^{opt})/\epsilon] \label{eq:opt-thm2} 
  \end{eqnarray}
  \end{thm} 
  
  Suppose our implementation of the stochastic shortest path problem every frame
 is accurate enough to ensure $\delta =0$. Then 
 from \eqref{eq:opt-thm2} and \eqref{eq:opt-thm1}, 
 the time average of $y_0(t)$  can be made arbitrarily close to (or below) 
 $y_0^{opt}$
 as $V$ is increased, with a tradeoff in average queue size that is linear in $V$. 
 The dependence on the $\phi$ parameter is also apparent:  While we desire $\phi$ to be small to minimize the disruptions
 due to forced renewals, a small value of $\phi$ implies a larger value of $B/\phi$ in  \eqref{eq:opt-thm2} and 
  \eqref{eq:opt-thm1}.  Note also that the average size of each stability-constrained queue affects its average delay, and the average size of
  each virtual queue affects the convergence time required for its constraint to be closely met. 
 
 \subsection{Proof of Theorem \ref{thm:1}} 
 
 We first prove \eqref{eq:opt-thm1}, and then \eqref{eq:opt-thm2}. 
 
 \begin{proof} (Theorem \ref{thm:1} part 1---Queue Bounds) 
 Let $t_r$ be the start of a renewal time. 
From \eqref{eq:delta-opt} and \eqref{eq:rhs-drift} we have: 
\begin{eqnarray}
\Delta(t_r) + V \expect{\sum_{\tau=t_r}^{t_r+T_r-1} y_0(\tau)\left|\right.\bv{\Theta}(t_r)} \leq
B/\phi^2 + C  \nonumber \\
+ \expect{D^*(\bv{\Theta}(t_r)) + V\sum_{\tau=t_r}^{t_r+T_r-1} y_0^*(\tau)\left|\right.\bv{\Theta}(t_r)} \nonumber \\
+ \delta\sum_{n\in\script{N}}Q_n(t_r) + \delta\sum_{l=1}^LX_l(t_r)  
+ V\delta  \label{eq:drift-opt1} 
\end{eqnarray}
where $D^*(\bv{\Theta}(t_r))$ and $y_l^*(\tau)$ are for any alternative policy $\alpha^*(t)$. Using the fact that $|y_0^*(\tau) - y_0(\tau)| \leq 2\beta$ for
all $\tau$, and $\expect{T_r} = 1/\phi$, we have: 
\begin{eqnarray}
\Delta(t_r) \leq
B/\phi^2 + C  + 2V\beta/\phi \nonumber \\
+ \expect{D^*(\bv{\Theta}(t_r)) \left|\right.\bv{\Theta}(t_r)} \nonumber \\
+ \delta\sum_{n\in\script{N}}Q_n(t_r) + \delta\sum_{l=1}^LX_l(t_r)  
+ V\delta  \label{eq:drift-opt2} 
\end{eqnarray}

Now consider the $(\Omega, z)$-only policy $\alpha_2^*(t)$ from Assumption 2, which makes decisions independent of $\bv{\Theta}(t_r)$ 
to yield (using the definition of $D(\bv{\Theta}(t_r))$ in \eqref{eq:D-def}): 
\begin{eqnarray*}
\expect{D^*(\bv{\Theta}(t_r))|\bv{\Theta}(t_r)} \leq \frac{-\epsilon}{\phi}\left[\sum_{n\in\script{N}} Q_n(t_r) + \sum_{l=1}^LX_l(t_r)\right] 
\end{eqnarray*}
Substituting the above into the right-hand-side of \eqref{eq:drift-opt2} gives: 
\begin{eqnarray}
\Delta(t_r) \leq
B/\phi^2 + C  + V(2\beta/\phi + \delta) \nonumber \\
+ (\delta - \epsilon/\phi)\left[\sum_{n\in\script{N}}Q_n(t_r) + \sum_{l=1}^LX_l(t_r)\right]    \label{eq:drift-opt3} 
\end{eqnarray}
Taking expectations of the above and using the definition of $\Delta(t_r)$ gives: 
\begin{eqnarray*}
\expect{L(t_{r+1})} - \expect{L(t_r)}  \leq
B/\phi^2 + C  + V(2\beta/\phi + \delta) \nonumber \\
+ (\delta - \epsilon/\phi)\left[\sum_{n\in\script{N}}\expect{Q_n(t_r)} + \sum_{l=1}^L\expect{X_l(t_r)}\right]  
\end{eqnarray*}
Summing the above over $r \in \{0, \ldots, R-1\}$ (for some positive integer $R$), dividing by $R$, 
and using the fact that $\expect{L(t_0)} = 0$ gives: 
\begin{eqnarray*}
\frac{\expect{L(t_{R})}}{R}  \leq
B/\phi^2 + C  + V(2\beta/\phi + \delta) \nonumber \\
+ \frac{(\delta - \epsilon/\phi)}{R}\sum_{r=0}^{R-1}\left[\sum_{n\in\script{N}}\expect{Q_n(t_r)} + \sum_{l=1}^L\expect{X_l(t_r)}\right]    
\end{eqnarray*}
Rearranging terms and using $\expect{L(t_R)}\geq 0$ and $\epsilon > \phi\delta$ proves \eqref{eq:opt-thm1}.  
While \eqref{eq:opt-thm1} samples only at the start of renewal frames, it can easily be used to show all queues are strongly stable (recall that the maximum queue change over any slot is bounded, and frame sizes are geometrically distributed with average $1/\phi$). Hence,
by stability theory in \cite{sno-text} we know all desired inequality constraints are met.  
\end{proof} 

\begin{proof} (Theorem \ref{thm:1} part 2 --- Performance Bound) 
Define probability $\theta \defequiv \delta\phi/\epsilon$. 
This is  a valid probability because $\epsilon>\phi\delta$ by assumption.  
We consider a new policy $\alpha^*(t)$ implemented 
over the frame $\tau \in \{t_r, \ldots, t_r+T_r-1\}$.  The policy $\alpha^*(t)$ is  a randomized
mixture of the $(\Omega,z)$-only policies from Assumptions 1 and 2:   
At the start of the frame, independently 
flip a biased coin with probabilities $\theta$ and $1-\theta$, and carry out one
of the two following policies for the full duration of the renewal interval:  
\begin{itemize} 
\item With probability $\theta$: Use policy $\alpha_2^*(t)$  from Assumption 2
for the duration of the renewal frame, 
which yields \eqref{eq:a2-1}-\eqref{eq:a2-2}.
\item With probability $1-\theta$: Use policy $\alpha_1^*(t)$  from Assumption 1 
for the duration of the renewal frame, 
which yields \eqref{eq:a1-1}-\eqref{eq:a1-3}. 
\end{itemize} 

Note that this policy $\alpha^*(t)$ is independent of $\bv{\Theta}(t_r)$. 
With $\alpha^*(t)$,  from \eqref{eq:a1-1} we have: 
\begin{eqnarray}
\expect{\sum_{\tau=t_r}^{t_r+T_r-1} y_0^*(\tau)|\bv{\Theta}(t_r)} \leq \frac{\theta \beta + (1-\theta)y_0^{opt}}{\phi}   \label{eq:part1} 
\end{eqnarray}
We also have from \eqref{eq:a1-2}-\eqref{eq:a1-3} and \eqref{eq:a2-1}-\eqref{eq:a2-2}: 
\begin{eqnarray}
\expect{\sum_{\tau=t_r}^{t_r+T_r-1}y_l^*(\tau)|\bv{\Theta}(t_r)} &\leq& \frac{-\theta\epsilon}{\phi} = -\delta \: \: \forall l \in \{1, \ldots, L\}  \nonumber \\
\expect{\sum_{\tau=t_r}^{t_r+T_r-1} d_n^*(\tau)|\bv{\Theta}(t_r)} &\leq& \frac{-\theta \epsilon}{\phi} = -\delta \: \:  \forall n \in \script{N}  \label{eq:part3} 
\end{eqnarray}
Plugging \eqref{eq:part1}-\eqref{eq:part3}  into \eqref{eq:drift-opt1} yields: 
\begin{eqnarray*}
\Delta(t_r) + V \expect{\sum_{\tau=t_r}^{t_r+T_r-1} y_0(\tau)\left|\right.\bv{\Theta}(t_r)} &\leq&
B/\phi^2 + C + V\delta  \\
&& \hspace{-.3in} + \frac{V}{\phi}[\theta \beta + (1-\theta)y_0^{opt}]
\end{eqnarray*}
Taking expectations gives: 
\begin{eqnarray*}
\expect{L(t_{r+1})} - \expect{L(t_r)}  + V \expect{\sum_{\tau=t_r}^{t_r+T_r-1} y_0(\tau)} \leq \\
B/\phi^2 + C + V\delta  + \frac{V}{\phi}[\theta \beta + (1-\theta)y_0^{opt}]
\end{eqnarray*}
Summing over $r \in \{0, \ldots, R-1\}$ and dividing by $VR/\phi$ gives the following for all $R>0$: 
\[ \frac{1}{R/\phi}\expect{\sum_{\tau=0}^{t_R-1}y_0(\tau)} \leq [\theta\beta + (1-\theta)y_0^{opt} + \delta\phi]  + \frac{B/\phi + C\phi}{V} \]
Using $\theta = \delta\phi/\epsilon$ shows the right-hand-side of the above inequality is the same as the right-hand-side of the 
desired inequality \eqref{eq:opt-thm2}. Finally, because $|y_0(\tau)| \leq \beta$ for all $\tau$, and $\{T_r\}_{r=0}^{\infty}$ are i.i.d. geometric
random variables with mean
$\expect{T_r} = 1/\phi$, it can be shown that (see Appendix A): 
\[ \limsup_{R\rightarrow\infty} \frac{\expect{\sum_{\tau=0}^{t_R-1}y_0(\tau)}}{R/\phi} \geq \limsup_{t\rightarrow\infty} \frac{1}{t}\sum_{\tau=0}^{t-1} \expect{y_0(\tau)} \]
 \end{proof}

 \section{Approximating the Stochastic Shortest Path Problem} \label{section:approx}    
 
 Consider now the stochastic shortest path problem \eqref{eq:dpp}.  Here we describe several approximation options. 
For simplicity, assume the state space $(\Omega(t) ,z(t))$ is finite, and the action space $\script{A}_{\Omega(t), z(t)}$ 
is finite for all $(\Omega(t), z(t))$. 
  Without loss of generality, assume we start at time $0$ and have (possibly non-zero) 
 backlogs $\bv{\Theta} = \bv{\Theta}(0)$.  Let $T$ 
 be the renewal interval size.  
 For every step $\tau \in \{0, \ldots, T-1\}$, define $c_{\bv{\Theta}}(\alpha(\tau), \Omega(\tau), z(\tau))$ as 
 the incurred cost assuming that the queue state at the beginning of the renewal is $\bv{\Theta}(0)$: 
\begin{eqnarray} 
c_{\bv{\Theta}}(\alpha(\tau), \Omega(\tau), z(\tau)) &\defequiv&  
  \sum_{n\in\script{N}} Q_n(0)\hat{d}_n(\alpha(\tau), \Omega(\tau), z(\tau)) \nonumber \\
 &&      +  \sum_{l=1}^L X_l(0)\hat{y}_l(\alpha(\tau), \Omega(\tau), z(\tau)) \nonumber \\
  &&     + V\hat{y}_0(\alpha(\tau), \Omega(\tau), z(\tau)) \label{eq:cdef} 
\end{eqnarray}

Let $\alpha^{ssp}(\tau)$ denote the optimal control action on slot $\tau$ for solving the stochastic shortest
path problem, given that the controller
 first observes $\Omega(\tau)$ and $z(\tau)$. 
 Define $\tilde{\script{Z}} \defequiv \script{Z} \cup \{renewal\}$, where we have added a new state ``$renewal$'' to represent
 the renewal state, which is the termination state of the stochastic shortest path problem.  Appropriately 
 adjust the transition probabilities $P_{ij}(\alpha, \Omega)$ to account 
 for this new state \cite{bertsekas-dp}\cite{bertsekas-neural}. 
 Define $\bv{J}  = (J_z)|_{z\in\tilde{\script{Z}}}$ 
 as a vector of optimal costs, where $J_z$ is the minimum
 expected sum cost to the renewal state given that we start in state $z$, and $J_{renewal}  = 0$. 
 By basic dynamic programming theory \cite{bertsekas-dp}\cite{bertsekas-neural}, the optimal 
 control action on each slot $\tau$ (given $\Omega(\tau)$ and $z(\tau)$)  
 is: 
\begin{eqnarray}
 & \hspace{-.6in} \alpha(\tau) =  \arg\min_{\alpha \in \script{A}_{\Omega(\tau), z(\tau)}  }    \left[ c_{\bv{\Theta}}(\alpha, \Omega(\tau), z(\tau)) + \right. \nonumber \\
 & \hspace{+1.4in}\sum_{y \in \tilde{\script{Z}}} P_{z(\tau), y}(\alpha, \Omega(\tau))J_y ]  \label{eq:I-choice} 
\end{eqnarray}
  
 This policy is easily implemented provided that the $J_z$ values are known.  It is well known that the $\bv{J}$ vector
 satisfies the following vector dynamic programming equation:\footnote{One can also derive \eqref{eq:vector-dp} 
 by defining a value function $H(z, \Omega)$, writing the Bellman equation in terms of $H(z(t+1), \Omega(t+1))$,  
 taking an expectation
 with respect to the  i.i.d. $\Omega(t)$, $\Omega(t+1)$, and 
 defining  $J(z) \defequiv \mathbb{E}_{\Omega(t)}\{H(z, \Omega(t))\}$.} 
 \begin{equation} \label{eq:vector-dp} 
 \bv{J} = \expect{\min_{\alpha_z\in \script{A}_{\Omega, z}}\left[ \bv{c}_{\bv{\Theta}}(\alpha_z, \Omega) + P(\alpha_z, \Omega)\bv{J} \right]  }  
 \end{equation} 
 where we have used an \emph{entry-wise min} (possibly with different $\alpha_z$ actions being used for minimizing each entry $z \in \tilde{\script{Z}}$).
 Further, $\bv{c}_{\bv{\Theta}}(\alpha_z, \Omega)$ is defined as a vector with 
 entries $c_{\bv{\Theta}}(\alpha_z, \Omega, z)$, and 
 $P(\alpha_z, \Omega) = (P_{zy}(\alpha_z, \Omega))$ is the matrix of transition probabilities
 for  $\Omega$ and control
 action $\alpha_z$.  The expectation in \eqref{eq:vector-dp} is over the distribution of the i.i.d. process $\Omega$. 
 Because $\Omega(t)$ has the structure $\Omega(t) = [\omega(t),  \phi(t)]$, 
 where $\omega(t)$ is the  random outcome for slot $t$ and $\phi(t)$ is an independent Bernoulli process that has
 forced renewals with probability $\phi$, we can re-write the above vector equation as: 
 \begin{eqnarray} 
 \bv{J} = \phi\expect{ \min_{\alpha_z \in \script{A}_{[\omega,1],z}} \bv{c}_{\bv{\Theta}}^{(1)}(\alpha_z, \omega)} +  \nonumber \\
 (1-\phi)\expect{ \min_{\alpha_z\in \script{A}_{[\omega,0],  z}}\left[ \bv{c}_{\bv{\Theta}}^{(0)}(\alpha_z, \omega) + P^{(0)}(\alpha_z, \omega)\bv{J} \right]} \label{eq:vector-dp2} 
 \end{eqnarray} 
 where: 
 \begin{eqnarray*}
 \bv{c}_{\bv{\Theta}}^{(1)}(\alpha_z, \omega) &\defequiv& \bv{c}_{\bv{\Theta}}(\alpha_z, [\omega, 1]) \\
 \bv{c}_{\bv{\Theta}}^{(0)}(\alpha_z, \omega) & \defequiv& \bv{c}_{\bv{\Theta}}(\alpha_z, [\omega, 0]) \\
 P^{(0)}(\alpha_z, \omega) &\defequiv& P(\alpha_z, [\omega, 0])
 \end{eqnarray*}
 
 We assume the transition probabilities $P^{(0)}(\alpha_z, \omega)$ are known (recall that these are indeed
 known binary values as described in the model of Section 
 \ref{section:notation}).    We next show how to compute an approximation of $\bv{J}$
 based on random samples of $\omega(t)$ and using a classic Robbins-Monro iteration.

 
 \subsection{Estimation Through Random i.i.d. Samples} \label{section:iidsamples} 
 
 Suppose we have an infinite sequence of random variables arranged in batches with batch size $W$, with 
 $\omega_{kw}$ denoting the $w$th  sample of batch $k$. 
 All random variables are i.i.d. with probability distribution the same as $\omega(t)$, and all are
 independent of the queue state $\bv{\Theta}$ that is used
 for this stochastic shortest path problem. 
 Consider the following two mappings $\Psi$ and $\tilde{\Psi}$ from a $\bv{J}$ vector to another $\bv{J}$ vector, 
 where the second is implemented with respect to a particular batch $k$: 
 \begin{eqnarray}
 && \hspace{-.2in}  \Psi \bv{J} \defequiv  \phi\expect{\min_{\alpha_z\in \script{A}_{[\omega, 1]}, z} \bv{c}_{\bv{\Theta}}^{(1)}(\alpha_z, \omega)} + \nonumber \\
&& \hspace{-.2in}  (1-\phi)\expect{\min_{\alpha_z \in \script{A}_{[\omega,0],z}}\left[ \bv{c}_{\bv{\Theta}}^{(0)}(\alpha_z, \omega) + P^{(0)}(\alpha_z, \omega)\bv{J}  \right]} \label{eq:psi} \\
&&  \hspace{-.2in} \tilde{\Psi} \bv{J} \defequiv  \phi\frac{1}{W}\sum_{w=1}^W \min_{\alpha_z \in \script{A}_{[\omega_{kw}, 1], z}}
\bv{c}_{\bv{\Theta}}^{(1)}(\alpha_z, \omega_{kw}) +  \nonumber \\
&& \hspace{-.2in} (1-\phi)\frac{1}{W}\sum_{w =1}^{W}\min_{\alpha_z \in \script{A}_{[\omega_{kw}, 0], z}}\left[ \bv{c}^{(0)}_{\bv{\Theta}}(\alpha_z, \omega_{kw}) +\right. \nonumber \\
&& \hspace{+1.5in} \left. P^{(0)}(\alpha_z, \omega_{kw})\bv{J}  \right] \label{eq:psi-tilde} 
 \end{eqnarray}
 where the min is  entrywise over each vector entry. The expectation in \eqref{eq:psi} is implicitly conditioned on
 a given $\bv{\Theta}$ vector, and is with respect to the random $\omega$, which is independent of 
 $\bv{\Theta}$.  We note that both $\Psi\bv{J}$ and $\tilde{\Psi}{\bv{J}}$ are vectors with size determined by the 
 size of the state space $\script{Z}$.  For a system with $K$ delay-constrained queues, the size of $\script{Z}$ is exponential 
 in $K$.  Thus, any computation of the map $\Psi{\bv{J}}$ or $\tilde{\Psi}{\bv{J}}$ must update a number of entries that is 
 exponential in $K$.  This is why we desire $K$ to be small, even though the number of stability-constrained queues 
 $N$ can be large.

 The mapping $\Psi$ cannot be implemented without knowledge of the distribution of $\omega$ (so that the expectation
 can be computed), whereas the mapping $\tilde{\Psi}$ can be implemented as a ``simulation'' over the $W$ random
 samples $\omega_{kw}$ (assuming such samples can be generated or obtained).  However, the expected
 value of $\tilde{\Psi} \bv{J}$ is exactly equal to $\Psi \bv{J}$.  Thus, given an initial vector $\bv{J}_k$ for use in step $k$, 
 we can  write
 $\tilde{\Psi}\bv{J}_k = \Psi\bv{J}_k + \bv{\eta}_k$,
 where $\bv{\eta}_k$ is a zero-mean vector random variable.  Specifically, 
 the vector $\bv{\eta}_k$ satisfies: 
 \[ \expect{\bv{\eta}_k\left|\right.\bv{J}_k}  = \bv{0} \]
 Thus, while the 
 vector $\bv{\eta}_k$ is \emph{not} independent
 of $\bv{J}_k$, each entry is \emph{uncorrelated} with any deterministic function of $\bv{J}_k$.  That is, for each 
 entry $i$ and any deterministic function $f(\cdot)$ we have
 via iterated expectations: 
 \begin{equation} \label{eq:uncorrelated} 
  \expect{\eta_k[i] f(\bv{J}_k)}  = \expect{f(\bv{J}_k)\expect{\eta_k[i]\left|\right.\bv{J}_k}} = 0 
  \end{equation} 
  
  For $k \in \{0, 1, 2, \ldots\}$ we have the iteration: 
\begin{equation} \label{eq:robbins-monro} 
 \bv{J}_{k+1} = \frac{1}{k+1}\tilde{\Psi}\bv{J}_k + \frac{k}{k+1}\bv{J}_k 
 \end{equation} 
 This iteration is a classic \emph{Robbins-Monro} stochastic approximation algorithm. It can be
 shown that the $\bv{J}$ vector remains deterministically bounded for all $k$ (see Appendix B), 
 and that  $\Psi$ and $\tilde{\Psi}$ satisfy the requirements of Proposition 4.6 in Section 4.3.4 of \cite{bertsekas-neural}. 
Thus the above iteration is in the standard form for stochastic approximation theory, and ensures that: 
\[ \lim_{k\rightarrow\infty} \bv{J}_k = \bv{J}^* \: \: \mbox{ with prob. 1}  \]
where $\bv{J}^*$ is the cost vector associated with the optimal stochastic shortest path problem, that is, it is the solution
to (\ref{eq:vector-dp2}) and thus satisfies $\bv{J}^* = \Psi \bv{J}^*$. This holds for any batch size $W$ (including the 
simplest case $W=1$), 
although taking larger 
batches reduces the variance of the per-batch estimation and may improve overall convergence speed. 

Unfortunately, the above does not specify how many iterations are needed to yield a close approximation 
to the $\bv{J}^*$ value.  The intuition is that we can run the iterations for a ``large enough'' time, and hope that
we have obtained a close enough approximation to yield $C$ and $\delta$ values that can be used in Theorem \ref{thm:1}. 
\footnote{We note that an earlier version of this technical report attempted to overcome this challenge by analyzing error bounds
for the following alternative recursion: 
\begin{equation*}
 \bv{J}_{k+1} = \gamma\tilde{\Psi}\bv{J}_k + (1-\gamma)\bv{J}_k 
 \end{equation*} 
 While this can be shown to be a contraction under the norm $\norm{\bv{J}} = \max_i |J_i|$, our results erroneously claimed
 (in Lemma 7 of the old technical report) that it was a contraction under the norm $\norm{\bv{J}} = \max_i \sqrt{\expect{J_i^2}}$, 
 where $\bv{J}$ is treated as a random variable (the two norms are identical if $\bv{J}$ is treated as a constant).   The error
 was in: (i) erroneously passing expectations through max[] in a step that was skipped, 
 and (ii) using a result that $\norm{P\bv{J}} \leq \norm{\bv{J}}$ (which 
 is correct for a fixed transition probability matrix $P$), where one actually would need 
 $\norm{P(\bv{J})\bv{J}} \leq \norm{\bv{J}}$, where $P(\bv{J})$ 
 can be a function of $\bv{J}$, and this latter inequality does not necessarily hold.}

  \subsection{Recursive Methods for $\Psi$} 

Contraction results for general stochastic shortest path problems are given in \cite{bertsekas-neural}. 
The following is a related result with a simpler form that holds because of our forced renewal
structure.  For a given vector $\bv{X}$, define $\norm{\bv{X}}$ as the maximum absolute value of $\bv{X}$: 
\[ \norm{\bv{X}} \defequiv \max_i |X_i| \]
It is not difficult to show that for any vector $\bv{X}$ and any probability matrix $P$ with rows that sum to 1, 
and with a number of columns equal to the size of $\bv{X}$, 
we have $\norm{P\bv{X}} \leq \norm{\bv{X}}$.

\begin{lem} \label{lem:contraction-easy1}For any vectors $\bv{X}$, $\bv{Y}$ of the same size as $\bv{J}^*$, we have: 
\[ \norm{\Psi \bv{X} - \Psi \bv{Y}} \leq (1-\phi) \norm{\bv{X} - \bv{Y}} \]
\end{lem} 

\begin{proof} 
Note that for all $\alpha_{z}$ and all  $\omega$ we have: 
\begin{eqnarray*}
&&\hspace{-.3in} \bv{c}_{\bv{\Theta}}^{(0)}(\alpha_z, \omega) + P^{(0)}(\alpha_z,\omega) \bv{X}  \\
&=& \bv{c}_{\bv{\Theta}}^{(0)}(\alpha_z, \omega) + P^{(0)}(\alpha_z,\omega) \bv{Y}  + 
P^{(0)}(\alpha_z,\omega)(\bv{X} - \bv{Y}) \\
&\leq& \bv{c}_{\bv{\Theta}}^{(0)}(\alpha_z, \omega) + P^{(0)}(\alpha_z,\omega) \bv{Y}  + 
\norm{P^{(0)}(\alpha_z,\omega)(\bv{X} - \bv{Y})} \\
&\leq& \bv{c}_{\bv{\Theta}}^{(0)}(\alpha_z, \omega) + P^{(0)}(\alpha_z,\omega) \bv{Y}  + \norm{\bv{X} - \bv{Y}} 
\end{eqnarray*}
Define vector $\bv{c}_1$ by: 
\[ \bv{c}_1 \defequiv \phi\expect{\min_{\alpha_z \in \script{A}_{[\omega,1], z}}\bv{c}_{\bv{\Theta}}^{(1)}(\alpha_z,\omega)} \]
Therefore: 
\begin{eqnarray*}
\Psi \bv{X} &=& \bv{c}_1 + (1-\phi) \times \\
&& \expect{\min_{\alpha_z\in\script{A}_{[\omega, 0], z}}\left[ \bv{c}_{\bv{\Theta}}^{(0)}(\alpha_z,\omega) + P^{(0)}(\alpha_z, \omega)\bv{X} \right]} \\
&\leq& c_1 + (1-\phi) \times \\
&& \expect{\min_{\alpha_z\in\script{A}_{[\omega, 0], z}}\left[ \bv{c}_{\bv{\Theta}}^{(0)}(\alpha_z,\omega) + P^{(0)}(\alpha_z, \omega)\bv{Y} \right]}  \\
&& + (1-\phi) \norm{\bv{X} - \bv{Y}} \\
&=& \Psi \bv{Y} + (1-\phi)\norm{\bv{X} - \bv{Y}} 
\end{eqnarray*}

By switching the roles of $\bv{X}$ and $\bv{Y}$ it can similarly be shown: 
\[ \Psi\bv{Y} \leq \Psi \bv{X} + (1-\phi)\norm{\bv{X} - \bv{Y}} \]
The result follows. 
\end{proof} 

This simple result yields the following approximation bounds for $k$ iterations of the map $\Psi$.
Define $\bv{J}_0$ as an initial guess of $\bv{J}^*$, and for $k \in \{1, 2, 3, \ldots\}$ define $\bv{J}_k = \Psi \bv{J}_{k-1}$. 

\begin{lem}\label{lem:contraction-easy2} For any initial vector $\bv{J}_0$ and any $k \in\{0, 1, 2,\ldots\}$ we have: 
\[ \norm{\bv{J}_k - \bv{J}^*} \leq (1-\phi)^k\norm{\bv{J}_0 -\bv{J}^*} \]
\end{lem}
\begin{proof} 
Recall that $\Psi \bv{J}_{k-1} = \bv{J}_k$ and $\Psi\bv{J}^* = \bv{J}^*$.  Then: 
\begin{eqnarray*}
\norm{\bv{J}_k - \bv{J}^*} &=& \norm{\Psi\bv{J}_{k-1} - \Psi\bv{J}^*} \\
&\leq& (1-\phi)\norm{\bv{J}_{k-1} - \bv{J}^*} 
\end{eqnarray*}
The result then follows easily by recursion. 
\end{proof}

Because the renewal frame size is independent of the policy, and has average $1/\phi$, 
it is not difficult to show that $\bv{J}^* \leq c_{max}/\phi$, where $c_{max}$ is the largest possible magnitude of 
of $c_{\bv{\Theta}}(\alpha(\tau), \Omega(\tau), z(\tau))$ for slot $\tau$ in the frame (such a constant
exists and is finite because of the boundedness assumptions). Therefore, defining $\bv{J}_0=\bv{0}$ and 
using Lemma \ref{lem:contraction-easy2} yields: 
\[ \norm{\bv{J}_k - \bv{J}^*} \leq (1-\phi)^kc_{max}/\phi \]
By the definition of $c_{\bv{\Theta}}(\cdot)$ in \eqref{eq:cdef}, it can be shown that $c_{max}$ is a sum of terms that are
proportional to $V$, $Q_n(t_r)$, and $Z_l(t_r)$.   Further, in 
Appendix C it is shown that the deviation in the optimal cost when \eqref{eq:I-choice} is used with an approximate
value $\bv{J}_k$, rather than $\bv{J}^*$, deviates from $\bv{J}^*$ by at most: 
\[ \frac{2(1-\phi)\norm{\bv{J}_k - \bv{J}^*}}{\phi}\]
Hence, the above two bounds can be used to compute a value $k$ that provides explicit approximation
values for $C$ and $\delta$ for use in Theorem \ref{thm:1}.

\subsection{Recursive Methods for $\tilde{\Psi}$}  

The difficulty in iterating the map $\Psi \bv{J}$ 
is that it requires full knowledge of the underlying 
probability distributions to compute the associated expectations.  
An approximation of this is to use $\tilde{\Psi}$ from (\ref{eq:psi-tilde}). Specifically, assume
we have $W$ i.i.d. samples $\omega_1, \ldots, \omega_W$.  Then the $\tilde{\Psi}$ function is: 
\begin{eqnarray*}
&&\hspace{-.2in}\tilde{\Psi} \bv{J} \defequiv  \phi\frac{1}{W}\sum_{w=1}^W \min_{\alpha_z \in \script{A}_{[\omega_{w}, 1], z}}
\bv{c}_{\bv{\Theta}}^{(1)}(\alpha_z, \omega_{w}) +  \nonumber \\
&&\hspace{-.2in} (1-\phi)\frac{1}{W}\sum_{w =1}^{W}\min_{\alpha_z \in \script{A}_{[\omega_{w}, 0], z}}\left[ \bv{c}^{(0)}_{\bv{\Theta}}(\alpha_z, \omega_{w}) + P^{(0)}(\alpha_z, \omega_{w})\bv{J}  \right] 
 \end{eqnarray*}

Define $\tilde{\bv{J}}_0$ as any initial vector, 
and for $k \in \{1, 2 ,3 , \ldots\}$  define $\tilde{\bv{J}}_k = \tilde{\Psi}\bv{J}_{k-1}$. 
Using the
same proof technique as Lemmas \ref{lem:contraction-easy1} and \ref{lem:contraction-easy2}, 
it is easy to show that for any $W>0$,  $\tilde{\Psi}$ is also a contraction that satisfies for any $\bv{X}$ and $\bv{Y}$: 
\[ \norm{\tilde{\Psi} \bv{X} - \tilde{\Psi} \bv{Y}} \leq (1-\phi)\norm{\bv{X} - \bv{Y}} \]
 Thus, it has a unique fixed point $\tilde{\bv{J}}^*$ satisfying 
$\tilde{\Psi}\tilde{\bv{J}}^* = \tilde{\bv{J}}^*$, and for all $k \in \{1, 2, 3, \ldots\}$ we have: 
\[ \norm{\tilde{\bv{J}}_k - \tilde{\bv{J}}^*} \leq (1-\phi)^k\norm{\tilde{\bv{J}}_0 - \tilde{\bv{J}}^*} \]

The value $\tilde{\bv{J}}^*$ is typically not the same as $\bv{J}^*$.  It represents the optimal cost vector in a modified
system where the $\omega$ vector is i.i.d. with the same distribution as the empirical average given over the $W$ samples. 
Intuitively, $\tilde{\bv{J}}^*$ becomes a better approximation for $\bv{J}^*$ when the number of samples $W$ is large.

 \subsection{Sampling From the Past and Delayed Queue Analysis} \label{subsection:delayed-queue}

It remains to be seen how one can obtain the required i.i.d. samples without knowing the probability distribution 
for $\omega$.  In this subsection,  we describe a technique that uses previous samples
 of the $\omega(\tau)$ values.

 We first obtain a collection of $W$ i.i.d. samples of $\omega(t)$.  Consider a given renewal time $t_r$,
 and suppose that the time $t_r$ is large enough so that we can obtain $W$ samples according to the 
 following procedure: Let $\omega_1 \defequiv \omega(t_r)$, 
 $\omega_2 \defequiv \omega(t_r-1)$, $\omega_3 \defequiv \omega(t_r-2), \ldots, \omega_W \defequiv \omega(t_r-W+1)$. 
Because $\omega(t)$ is i.i.d. over slots (and because our renewal times are chosen randomly and independently),
it is easy to see that $\{\omega_1, \ldots, \omega_W\}$ form an i.i.d. sequence.  

A subtlety now arises:  Even though the $\{\omega_1, \ldots, \omega_W\}$ sequence is i.i.d., these samples
are \emph{not} independent of the queue backlog $\bv{\Theta}(t_r)$ at the beginning of the renewal.  This is 
because these values have influenced the queue states.  This makes it challenging to directly implement 
a Robbins-Monro iteration.  Indeed, the
expectation in \eqref{eq:psi} can be viewed as a conditional expectation given a certain queue backlog
at the beginning of the renewal interval, which is $\bv{\Theta}(t_r)$ for the $r$th renewal. This conditioning
does not affect \eqref{eq:psi} when $\omega(t)$ is chosen independently of initial queue backlog, and so  the random
samples in \eqref{eq:psi-tilde} are also assumed to be chosen independent of the initial queue backlog, 
which is not the case if we sample from the past. 

 To avoid this 
difficulty and ensure the samples are both i.i.d. and independent
of the queue states that form the weights in our stochastic shortest path problem, we use a \emph{delayed queue
analysis} as in \cite{neely-mwl-arxiv}.
Let $t_{start}$ denote the slot on which sample $\omega_W$ is taken, and let $\bv{\Theta}(t_{start})$ represent the 
queue backlogs at that time.  It follows that the i.i.d. samples are also independent of $\bv{\Theta}(t_{start})$.  
Hence, the bounds derived for the iteration technique in the previous section can be applied when the 
iterates use $\bv{\Theta}(t_{start})$ as the backlog vector.  Let $\bv{J}_{\bv{\Theta}(t_r)}$ denote the optimal 
solution to the problem \eqref{eq:vector-dp} for a queue backlog $\bv{\Theta}(t_r)$ at the beginning 
of our renewal time $t_r$, and let $\bv{J}_{\bv{\Theta}(t_{start})}$ 
denote the corresponding optimal solution for a problem that starts with initial 
queue backlog $\bv{\Theta}(t_{start})$.   Then there are $W-1$ slots in between $t_{start}$ and $t_r$. 
Because the 
maximum change in any queue on one slot is  bounded by $\beta$, we want to claim that an algorithm which 
computes the stochastic shortest path using the $\bv{\Theta}(t_{start})$ queue values gives a result that is 
within an additive constant of the algorithm which uses $\bv{\Theta}(t_r)$. Such an additive constant can be viewed
as the $C$ constant in Theorem \ref{thm:1}. 
 This can be justified using the next lemma,
which bounds the deviation of the optimal costs associated with two general queue backlog vectors.

Let $\bv{\Theta}_1$ and $\bv{\Theta}_2$ be two different queue backlog vectors, and let 
$\bv{J}_{\bv{\Theta}_1}$ and $\bv{J}_{\bv{\Theta}_2}$ represent the optimal frame costs corresponding to $\bv{\Theta}_1$ and 
$\bv{\Theta}_2$, respectively.  Define the constant $\theta$ as follows: 
\begin{equation} \label{eq:beta-theta} 
\theta \defequiv \sup_{\alpha_z, \Omega} \norm{\bv{c}_{\bv{\Theta}_1}(\alpha_z, \Omega) - \bv{c}_{\bv{\Theta}_2}(\alpha_z, \Omega)} 
\end{equation} 
where $\bv{c}_{\bv{\Theta}}(\alpha_z, \Omega)$ is the vector, indexed by $z$, with the $z$th entry  given by \eqref{eq:cdef} using
backlog vector $\bv{\Theta}$.  Note from \eqref{eq:cdef} that $\theta$ is independent of $V$ (as the $V$ term in 
\eqref{eq:cdef} cancels out in the subtraction), and is proportional to the maximum penalty value times the 
maximum \emph{difference} in any queue backlog entry in $\bv{\Theta}_1$ and its corresponding entry in 
$\bv{\Theta}_2$.  Thus $\theta$ is also independent of the actual size of the backlog vectors, and depends only
on their \emph{difference}, being proportional to $W\beta$.

\begin{lem}  \label{lem:j-diff} For the vectors $\bv{\Theta}_1$ and $\bv{\Theta}_2$, and for the $\theta$ value
defined in \eqref{eq:beta-theta}, we have: 

(a) The difference between $\bv{J}_{\bv{\Theta}_1}$ and $\bv{J}_{\bv{\Theta}_2}$ satisfies: 
\[   \norm{\bv{J}_{\bv{\Theta}_1}  - \bv{J}_{\bv{\Theta}_2}} \leq \frac{\theta}{\phi}  \]

(b) Let $\alpha_1(t)$ denote the policy decisions at time $t$ under the policy that makes optimal decisions subject
to queue backlogs $\bv{\Theta}_1$, and define $\bv{J}_{21}^{mis}$ as the expected sum cost over a frame of a 
\emph{mismatched policy} that incurs costs according to backlog vector $\bv{\Theta}_2$ but makes decisions
according to $\alpha_1(t)$ (and hence has the same decisions as the optimal policy for $\bv{\Theta}_1$). 
Then: 
\[ \bv{J}_{\bv{\Theta}_2} \leq \bv{J}_{21}^{mis} \leq \bv{J}_{\bv{\Theta}_1} + \bv{1}\frac{\beta}{\phi} \]
where $\bv{1}$ is a vector of all $1$ values with the same dimension as $\bv{J}_{\bv{\Theta}_1}$. 
\end{lem}   
\begin{proof}
 Omitted for brevity (see Appendix D).   
\end{proof} 

\section{Simulation} \label{section:simulation} 

\begin{figure}[b]
  \centering
  \includegraphics{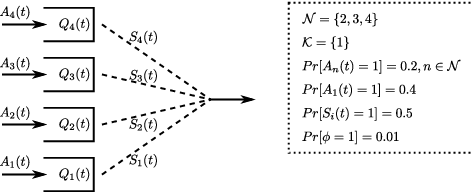}
  \caption{A network with 1 delay-constrained queue (queue 1), and 3 stability-constrained queues.}
  \label{fig:simplenetwork}
\end{figure}

In this section, we simulate the frame-based drift-plus-penalty algorithm in Section \ref{section:dpp-alg} 
for the simple network in Fig. \ref{fig:simplenetwork}. The algorithm utilizes the classic Robbins-Monro iteration, based on samples from
the past,  to approximate 
the weighted stochastic shortest path problem \eqref{eq:vector-dp2}. This is because solving 
\eqref{eq:vector-dp2} exactly is computationally expensive, would require full probability knowledge, and 
may not be practical for implementation.

The network in Fig. \ref{fig:simplenetwork} consists of one delay-constrained queue and three stability-constrained queues, so that $\mathcal{K} = \{ 1 \}$ and $\mathcal{N} = \{ 2, 3, 4 \}$. The size of the delay-constrained queue is limited to $b = 10$ packets. Random packet arrivals are i.i.d. Bernoulli processes with $Pr[ A_n(t) = 1 ] = 0.2$ for $n \in \mathcal{N}$ and $Pr[ A_1(t) = 1 ] = 0.4$. Each network channel is a binary state and is active (ON-state) with probability $Pr[ S_i(t) = 1 ] = 0.5$ for $i \in \mathcal{N} \cup \mathcal{K}$. The force renewal probability is $Pr[\phi(t) = 1] = 0.01$.

In this simulation, we consider a problem of minimizing the average number of dropped packets. For the delay-constrained queue $Q_1(t)$, the average backlog is limited to 1.5. Define $y_0(t) = D_1(t)$ and $y_1(t) = Q_1(t) - 1.5$. Then an optimization for this simulation is
\begin{eqnarray*}
  \text{minimize}   & \bar{y}_0 \\
  \text{subject to} & \bar{y}_1 \leq 1.5 \\
                    & \bar{Q}_n < \infty ~ \text{for all} ~ n \in \{2, 3, 4 \}.
\end{eqnarray*}

The simulation follows the frame-based drift-plus-penalty algorithm in Section \ref{section:dpp-alg} 
with the Robbins-Monro iteration \eqref{eq:robbins-monro}. A batch size is set to be $W = 50$, so that we store the most recent $50$ samples (using 
less than 50 in the initial slots $\tau < 50$). Note that the number of samples is
 half of the average frame size, $1/\phi = 100$. Every forced renewal slot $t_r$, the algorithm uses the batch to approximate the 
 mapping $\tilde{\Psi} \bv{J}$ in \eqref{eq:psi-tilde}, and then updates $\bv{J}$ according to \eqref{eq:robbins-monro}. 
 After updating $\bv{J}$, every decision in frame $r$ is decided from the simple rule \eqref{eq:I-choice}. 
 Then all delay-constrained, stability-constrained, and virtual queues are updated as in 
 \eqref{eq:q-update}, \eqref{eq:z-update}, and \eqref{eq:x-update}.

\begin{figure}
  \centering
  \includegraphics[scale=0.27]{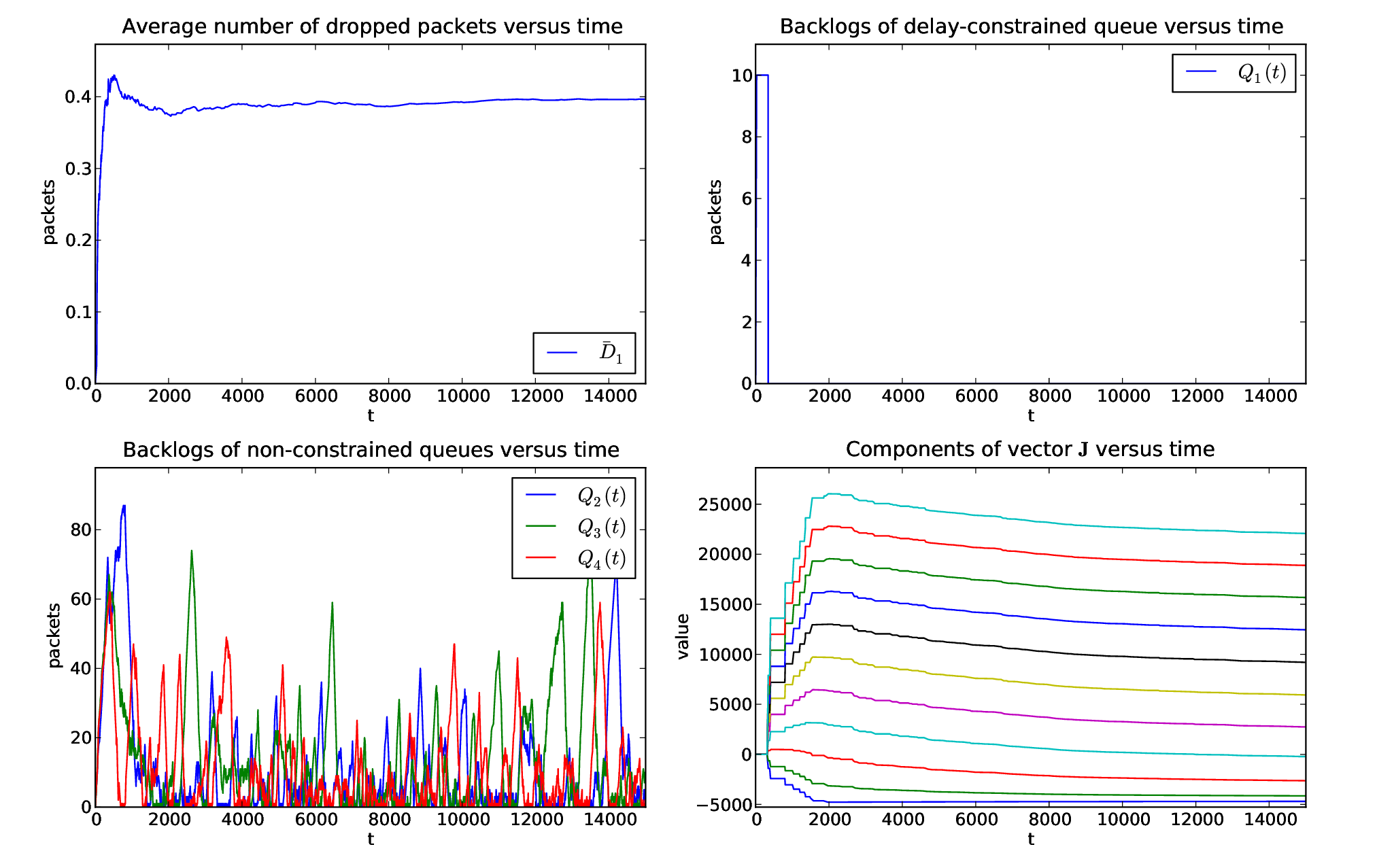}
  \caption{Behavior in the system with $V = 0$}
  \label{fig:finalV0}
\end{figure}

For a simple initial comparison, we use $V=0$, so the algorithm puts no weight on minimizing $\overline{y}_0$ and only attempts
to satisfy the desired constraints. Results from the algorithm until $1.5 \times 10^4$ slots are shown in Fig. \ref{fig:finalV0}. 
The system drops almost all packets in the delay-constrained queue (as expected),
making its average queue size approach zero, as shown in the top graphs of Fig. \ref{fig:finalV0}. 
All stability-constrained queues are stable, which is shown in the bottom-right of Fig. \ref{fig:finalV0}. The bottom-left of Fig. \ref{fig:finalV0} shows the convergence of $\bv{J}$. This illustrates that the algorithm yields a feasible solution.

\begin{figure}
  \centering
  \includegraphics[scale=0.27]{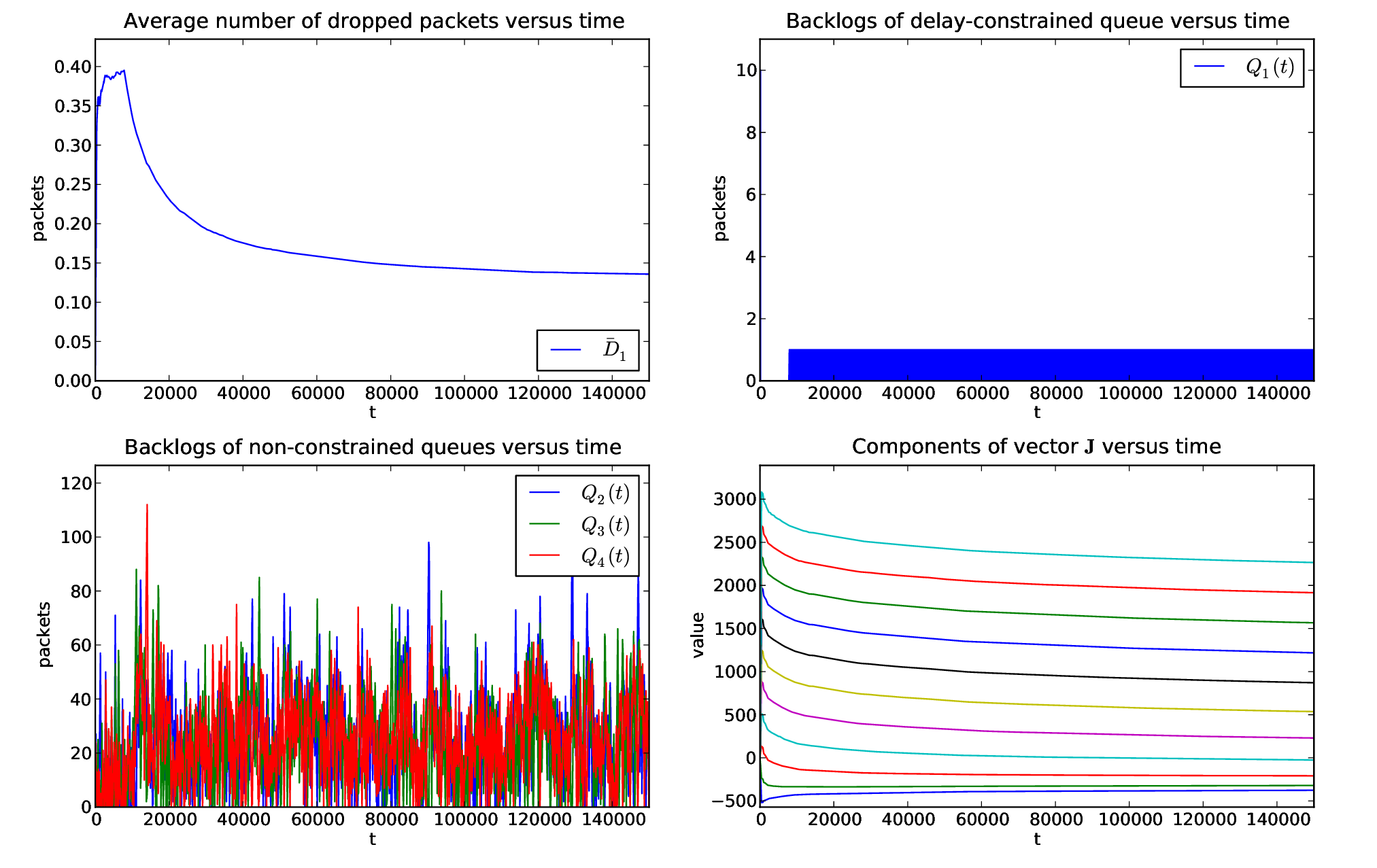}
  \caption{Behaviors of the system with $V = 100$}
  \label{fig:finalV100}
\end{figure}

\begin{figure}
  \centering
  \includegraphics[scale=0.37]{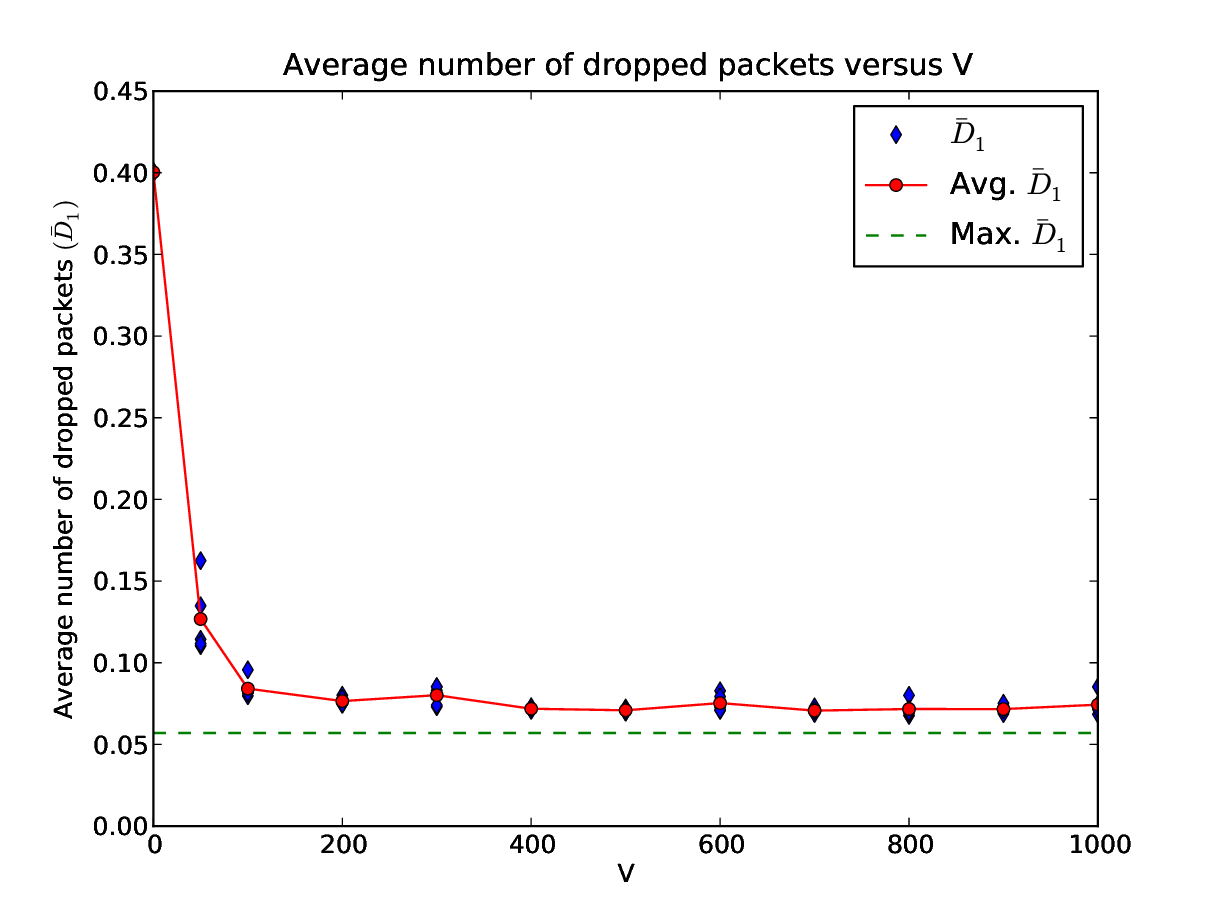}
  \caption{Average number of dropped packets versus $V$}
  \label{fig:performance}
\end{figure}

\begin{figure}
  \centering
  \includegraphics[scale=0.37]{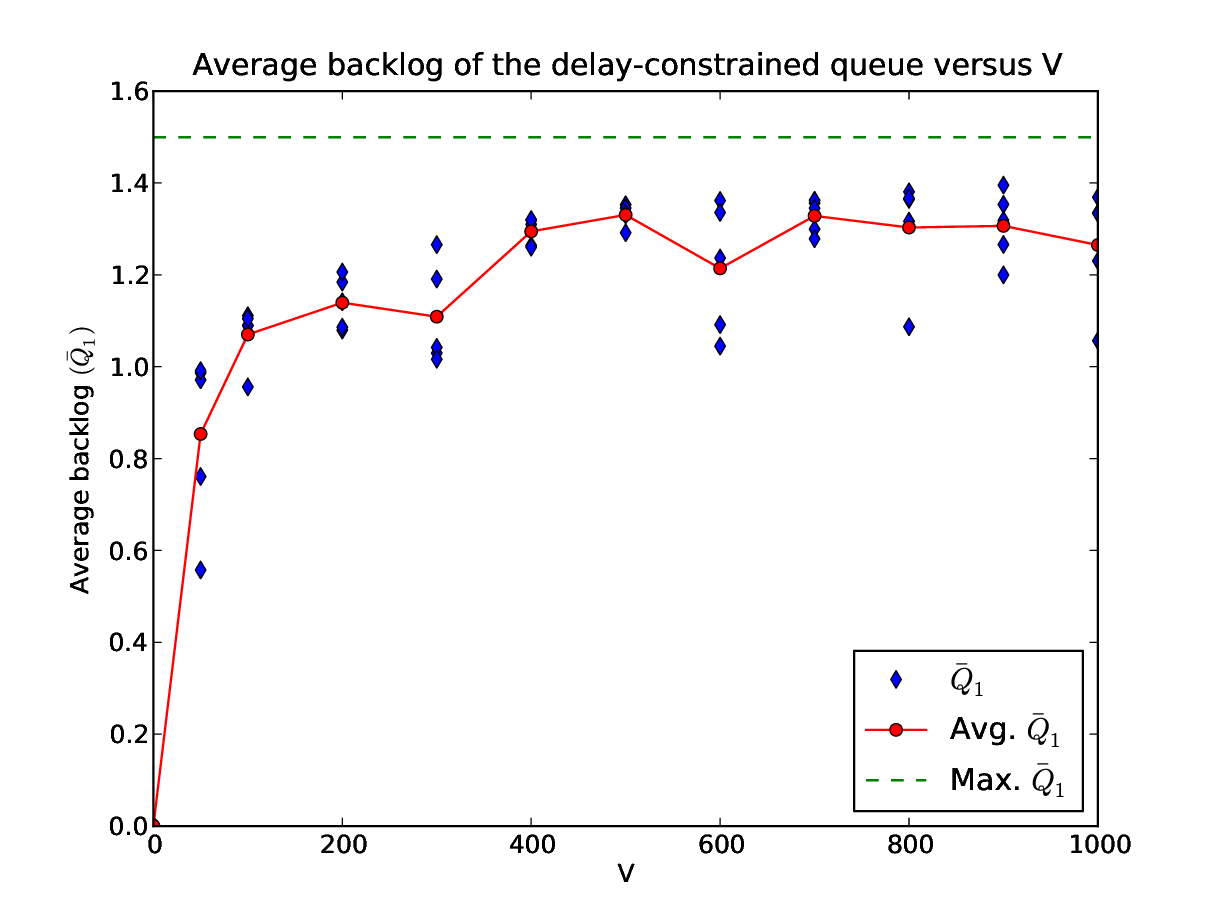}
  \caption{Average backlog of the delay-constrained queue versus $V$}
  \label{fig:backlog}
\end{figure}

We next use $V=100$, so the algorithm attempts to minimize dropping in queue 1.  Behaviors in the system for 
the first $1.5 \times 10^5$ slots are shown in Fig. \ref{fig:finalV100}. The figure shows the convergence of the algorithm. After $10^6$ slots, the average number of dropped packets is $0.096$ and the average backlog of the delay-constrained queue is $0.956$. These values correspond to the
data points plotted for $V = 100$ in Figs. \ref{fig:performance} and \ref{fig:backlog}. 
Compared to the result from $V = 0$, the average number of dropped packets decreases, while the backlog increases as a result of more aggressive admission. In addition, the algorithm with $V = 100$ takes more time slots to converge.

Finally, the system is simulated for $V$ in the range from 
$0$ to $1000$, as shown in Figs. \ref{fig:performance} and \ref{fig:backlog}. Each value of $V$ is simulated over 5 independent runs. 
As $V$ is increased, we expect the average drop rate to converge to optimality, with a corresponding increase in average queue sizes for the 
stability-constrained queues.  This is exactly what happens.  After $10^6$ slots, the average number of dropped packets and the average number of backlogs are recorded. Then the average of the five values for each $V$ is calculated. Also, the nearest optimal solution (when $V = 10^4$) that we obtained is represented by dashed lines in both figures. In this case, the average number of dropped packets is $0.057$, and the average number of backlog is $\overline{Q}_1 = 1.499$. Note that in this case, the 
average queue size constraint $\overline{Q}_1 \leq 1.5$ is met with near equality, which is why the number of dropped
packets can be pushed down so far. Fig. \ref{fig:performance} illustrates the performance of the algorithm as $V$ varies. 

\section{Conclusions} 

We have developed an approach to the Markov Decision problems associated with a small number
 $K$ of delay-constrained wireless users and a (possibly large) number $N$ of stability-constrained queues.  
 Our formulation allows 
 optimization 
 of general penalty functions subject to general penalty constraints, such as minimizing average packet drops subject
 to average backlog and/or average delay constraints at the delay-constrained queues, and subject to 
 stability at the stability-constrained queues.  Our approach uses 
 a reduction to an online (unconstrained) weighted stochastic shortest path problem implemented over variable length frames.  
 This generalizes the class of max-weight network control policies to networks with Markov decisions. 
 The solution to the underlying stochastic shortest path problem has complexity that is exponential in the number
 of delay-constrained queues $K$, but polynomial in the number of delay-unconstrained queues $N$.  
 A Robbins-Monro approximation technique was used to develop several approximation algorithms for the stochastic
 shortest path problem. 
  The solution technique is general and
 extends to other network problems with stochastic decisions.

\section*{Appendix A} 

Here we show that if $\{T_r\}_{r=0}^{\infty}$ are i.i.d. geometric random variables with mean $\expect{T_r} = 1/\phi$, 
and if $|y_0(\tau)| \leq \beta$ for all $\tau$ (for some positive constant $\beta$), 
then: 
\begin{equation} \label{eq:appa} 
 \limsup_{R\rightarrow\infty} \frac{\expect{\sum_{\tau=0}^{t_R-1}y_0(\tau)}}{R/\phi} \geq \limsup_{t\rightarrow\infty} \frac{1}{t}\sum_{\tau=0}^{t-1} \expect{y_0(\tau)} 
 \end{equation} 
This  is used at the end of the proof of Theorem \ref{thm:1}. 

We have by the 
Law of Large Numbers: 
\[ \frac{t_R}{R} =  \frac{\sum_{r=0}^{R-1} T_r}{R} \rightarrow 1/\phi \: \: \mbox{with prob. 1} \]

First consider the case when $y_0(\tau) \geq 0$ with probability 1 for all $\tau$,  so that $0 \leq y_0(\tau) \leq \beta$ for all $\tau$. 
Let $R(t)$ be the number of renewal events that have occurred up to time $t$ (not counting the renewal 
at time $0$).  Then $R(t)/t \rightarrow \phi$ with probability $1$.  Fix a value $\epsilon>0$ such that $0 < \epsilon < \phi$. Define the following event $\chi(t)$: 
\[ \chi(t) \defequiv \left\{   \frac{R(t)}{t} <  \phi + \epsilon\right\} \]
Define $\chi^c(t)$ as the opposite event. 
Then $Pr[\chi^c(t)] \rightarrow 0$ as $t \rightarrow \infty$.    If $\chi(t)$ is true, then
$R(t) <  \lceil (\phi+\epsilon)t\rceil$ and so: 
\[ t <  t_{\lceil(\phi+\epsilon)t\rceil}  \: \mbox{ whenever $\chi(t)$ is true} \]
where we recall that $t_r$ is the time of the $r$th renewal event. 
Now for any time $t$ we have: 
\begin{eqnarray*}
 \frac{1}{t}\sum_{\tau=0}^{t-1} \expect{y_0(\tau)} 
&=& \expect{\frac{1}{t}\sum_{\tau=0}^{t-1} y_0(\tau) \left|\right. \chi(t)}Pr[\chi(t)] \\
&& + \expect{\frac{1}{t}\sum_{\tau=0}^{t-1} y_0(\tau)\left|\right.\chi^c(t)}Pr[\chi^c(t)]  \\
&\leq& \expect{\frac{1}{t}\sum_{\tau=0}^{t_{\lceil(\phi+\epsilon)t\rceil} - 1 } y_0(\tau) \left|\right.\chi(t)}Pr[\chi(t)] \\
&& +\beta Pr[\chi^c(t)] \\
&\leq& \frac{1}{t}\expect{\sum_{\tau=0}^{t_{\lceil(\phi+\epsilon)t\rceil}-1} y_0(\tau)} \\
&&+ \beta Pr[\chi^c(t)]
\end{eqnarray*}
where the final inequality holds because we have added the non-negative term: 
\[ \frac{1}{t}\expect{\sum_{\tau=0}^{t_{\lceil(\phi+\epsilon)t\rceil}-1} y_0(\tau)\left|\right.\chi^c(t)}Pr[\chi^c(t)] \]
Therefore: 
\begin{eqnarray*}
&&\hspace{-.4in} \frac{1}{t}\sum_{\tau=0}^{t-1} \expect{y_0(\tau)} \\
 &\leq& \frac{\lceil(\phi+\epsilon)t\rceil} {t}\frac{1}{\lceil(\phi+\epsilon)t\rceil}\expect{\sum_{\tau=0}^{t_{\lceil(\phi+\epsilon)t\rceil}-1} y_0(\tau)}  \\
 &&+ \beta Pr[\chi^c(t)]
\end{eqnarray*}
Taking limits yields: 
\begin{eqnarray*}
\limsup_{t\rightarrow\infty} \frac{1}{t}\sum_{\tau=0}^{t-1} \expect{y_0(\tau)} \leq (\phi+\epsilon)\limsup_{R\rightarrow\infty}\frac{1}{R}\expect{\sum_{\tau=0}^{t_R - 1} y_0(\tau)}
\end{eqnarray*}
The above holds for all $\epsilon$ such that $0 < \epsilon < \phi$.  Taking a limit as $\epsilon\rightarrow 0$ yields: 
\begin{eqnarray*}
\limsup_{t\rightarrow\infty} \frac{1}{t}\sum_{\tau=0}^{t-1} \expect{x_0(\tau)} \leq \phi\limsup_{R\rightarrow\infty}\frac{1}{R}\expect{\sum_{\tau=0}^{t_R- 1} y_0(\tau)}
\end{eqnarray*}
The reverse inequality can be proven similarly.  This establishes (\ref{eq:appa}) for the case when $y_0(\tau)$
is a non-negative process.  

For the case $|y_0(\tau)| \leq \beta$, but can take possibly negative values, 
we can define $\tilde{y}_0(\tau) \defequiv y_0(\tau) + \beta$.  Then 
we have $0 \leq \tilde{y}_0(\tau) \leq 2\beta$ for all $\tau$.  It follows that: 
\begin{eqnarray*}
\limsup_{t\rightarrow\infty} \frac{1}{t}\sum_{\tau=0}^{t-1} \expect{\tilde{y}_0(\tau)} \leq \phi\limsup_{R\rightarrow\infty}\frac{1}{R}\expect{\sum_{\tau=0}^{t_R- 1} \tilde{y}_0(\tau)}
\end{eqnarray*}
Subtracting $\beta$ from both sides of the above equality yields the result of (\ref{eq:appa}).

\section*{Appendix B} 

Here we show that if we use iteration \eqref{eq:robbins-monro} starting with any initial vector $\bv{J}_0$, then 
the norms $\norm{\bv{J}_k}$ of all iterates $\bv{J}_k$ are bounded, where we use the max-absolute value norm: 
\[ \norm{\bv{X}} \defequiv \max_i|X_i| \]

Consider the iteration: 
\[ \bv{J}_{k+1} = \frac{1}{k+1}\tilde{\Psi}\bv{J}_k + \frac{k}{k+1} \bv{J}_k \]
where:
\[ \tilde{\Psi} \bv{J}_k = \Psi \bv{J}_k + \bv{\eta}_k \]
where $\Psi$ is the map of (\ref{eq:psi}), and 
$\{\bv{\eta}_k\}_{k=0}^{\infty}$ is a sequence of  zero mean vector random variables, 
where each entry of $\bv{\eta}_k$ is uncorrelated with any deterministic function of 
$\bv{J}_k$.  We show that $\norm{\bv{J}_k}$ and $\norm{\bv{\eta}_k}$ are deterministically
bounded.   Define $c_{max}$ as the maximum absolute value of any term of
the $c_{\bv{\Theta}}^{(0)}(\alpha, \omega)$ or $c_{\bv{\Theta}}^{(1)}(\alpha,\omega)$ functions, under any $\alpha, \omega$.
This maximum is finite by the boundedness assumptions. 
Define $J_{max}\defequiv c_{max}/\phi$.  We claim that if 
If $\norm{\bv{J}_0} \leq J_{max}$, then: 
\[ \norm{\bv{J}_k} \leq J_{max} \: \: \: , \: \: \: \norm{\bv{\eta}_k} \leq 2J_{max} \]

\begin{proof} 
Suppose that $\norm{\bv{J}_k} \leq J_{max}$ for some iteration $k\geq0$ (it
holds by assumption for $k=0$).  We show that it also holds for $k+1$.  
By the update equations (\ref{eq:psi}) and (\ref{eq:psi-tilde}),  it is not difficult to show that: 
 \[ \max[\norm{\tilde{\Psi} \bv{J}_k} , \norm{\Psi\bv{J}_k}] \leq c_{max} + (1-\phi)J_{max} = J_{max}\]
Thus:  
\begin{eqnarray*}
 \norm{\bv{J}_{k+1}} &\leq& \frac{1}{k+1}\norm{\tilde{\Psi}\bv{J}_k} + \frac{k}{k+1}\norm{\bv{J}_k} \\
 &\leq&\frac{J_{max}}{k+1} + \frac{k}{k+1}J_{max}\\
 &=&  J_{max}  
 \end{eqnarray*}
 This proves the first part. To prove the second part, we have: 
 \begin{eqnarray*}
 \norm{\bv{\eta}_k} &=& \norm{\tilde{\Psi}\bv{J}_k - \Psi \bv{J}_k} \\
 &\leq& \norm{\tilde{\Psi}\bv{J}_k} + \norm{\Psi\bv{J}_k} \\
 &\leq&  2J_{max}
 \end{eqnarray*}
\end{proof} 

\section*{Appendix C}

We now show that an implementation that
chooses $\alpha(t)$ over a frame according to  (\ref{eq:I-choice}), using the $\bv{J}_k$ estimate instead of the 
optimal $\bv{J}^*$ vector, results in an approximation to the stochastic shortest path problem
that deviates by an amount that depends on $\norm{\bv{J}_k-\bv{J}^*}$. 

Claim: Suppose we choose $\alpha(t)$ according to (\ref{eq:I-choice}) over the course of a frame, 
using a vector $\bv{J}$ rather than $\bv{J}^*$. Let $\tilde{\bv{J}}(\bv{J})$ represent the expected sum cost over
the frame (given $\bv{J}$).  Then:
\begin{equation} \label{eq:cost-diff-app} 
   \norm{\tilde{\bv{J}}(\bv{J}) - \bv{J}^*} \leq \frac{2(1-\phi)\norm{\bv{J} - \bv{J}^*}}{\phi}   
   \end{equation}

\begin{proof} 
Let $\alpha(t)$ represent the control decision on slot $t$ made using the $\bv{J}$ vector, and let $\alpha^*(t)$ represent
the decision that would be made under the $\bv{J}^*$ vector. Then: 
\begin{eqnarray*}
&& \hspace{-.2in} \tilde{\bv{J}}(\bv{J}) 
= \phi \expect{\min_{\alpha_z\in \script{A}_{[\omega(t), 1], z}} \bv{c}_{\bv{\Theta}}^{(1)}(\alpha_z, \omega(t))} \\
&& + (1-\phi)\expect{\bv{c}_{\bv{\Theta}}^{(0)}(\alpha(t), \omega(t)) + P^{(0)}(\alpha(t), \omega(t))\tilde{\bv{J}}(\bv{J})\left|\right.\bv{J}}
\end{eqnarray*}
where the expectation is with respect to the random $\omega(t)$ outcome. Thus: 
\begin{eqnarray}
&&\hspace{-.3in}\tilde{\bv{J}}(\bv{J}) = \phi \expect{\min_{\alpha_z \in \script{A}_{[\omega(t), 1], \bv{z}}} \bv{c}_{\bv{\Theta}}^{(1)}(\alpha_z, \omega(t))} \nonumber \\
&& + (1-\phi)\expect{\bv{c}_{\bv{\Theta}}^{(0)}(\alpha(t), \omega(t)) + P^{(0)}(\alpha(t), \omega(t))\bv{J}\left|\right.\bv{J}} \nonumber \\
&&  + (1-\phi)\expect{P^{(0)}(\alpha(t), \omega(t))}(\tilde{\bv{J}}(\bv{J}) - \bv{J})  \label{eq:mango} 
\end{eqnarray}

Because $\alpha(t)$ minimizes the second term of the above equality, we have: 
\begin{eqnarray*}
&&\hspace{-.2in}(1-\phi)\expect{\bv{c}_{\bv{\Theta}}^{(0)}(\alpha(t), \omega(t)) + P^{(0)}(\alpha(t), \omega(t))\bv{J}\left|\right.\bv{J}}  \\
&&\leq (1-\phi)\expect{\bv{c}_{\bv{\Theta}}^{(0)}(\alpha^*(t), \omega(t)) + P^{(0)}(\alpha^*(t), \omega(t))\bv{J}\left|\right.\bv{J}} \\
&&= (1-\phi)\expect{\bv{c}_{\bv{\Theta}}^{(0)}(\alpha^*(t), \omega(t)) + P^{(0)}(\alpha^*(t), \omega(t))\bv{J}^*\left|\right.\bv{J}} \\
&& \hspace{+.1in}(1-\phi)\expect{P^{(0)}(\alpha^*(t), \omega(t))(\bv{J} - \bv{J}^*) \left|\right.\bv{J}}
\end{eqnarray*}
Combining the above with (\ref{eq:mango}) yields: 
\begin{eqnarray*}
&& \tilde{\bv{J}}(\bv{J}) \leq \bv{J}^* 
+ (1-\phi)\expect{P^{(0)}(\alpha(t), \omega(t))}(\tilde{\bv{J}}(\bv{J}) - \bv{J}) \\
&&+  (1-\phi)\expect{P^{(0)}(\alpha^*(t), \omega(t))}(\bv{J} - \bv{J}^*)
\end{eqnarray*}
However, we also know that $\bv{J}^* \leq \tilde{\bv{J}}(\bv{J})$. Therefore, 
using the fact that the expectation of a transition matrix is also a transition matrix, 
and that $\norm{P \bv{X}} \leq \norm{\bv{X}}$: 
\begin{eqnarray*}
  \norm{\tilde{\bv{J}}(\bv{J}) - \bv{J}^*} \leq (1-\phi) \norm{\bv{J} - \bv{J}^*} +   
  (1-\phi) \norm{\tilde{\bv{J}}(\bv{J}) - \bv{J}} \\
  \leq (1-\phi)[\norm{\bv{J}-\bv{J}^*} + \norm{\tilde{\bv{J}}(\bv{J}) - \bv{J^*}} + \norm{\bv{J}^*-\bv{J}}]
 \end{eqnarray*}
 Rearranging terms yields (\ref{eq:cost-diff-app}). 
  \end{proof}
 
 \section*{Appendix D} 

Here we prove Lemma \ref{lem:j-diff} of Section \ref{subsection:delayed-queue}, restated below for convenience:
For the vectors $\bv{\Theta}_1$ and $\bv{\Theta}_2$, and for the $\theta$ value
defined in \eqref{eq:beta-theta}, we have: 

(a) The difference between $\bv{J}_{\bv{\Theta}_1}$ and $\bv{J}_{\bv{\Theta}_2}$ satisfies: 
\[   \norm{\bv{J}_{\bv{\Theta}_1}  - \bv{J}_{\bv{\Theta}_2}} \leq \frac{\theta}{\phi}  \]

(b) Let $\alpha_1(t)$ denote the policy decisions at time $t$ under the policy that makes optimal decisions subject
to queue backlogs $\bv{\Theta}_1$, and define $\bv{J}_{21}^{mis}$ as the expected sum cost over a frame of a 
\emph{mismatched policy} that incurs costs according to backlog vector $\bv{\Theta}_2$ but makes decisions
according to $\alpha_1(t)$ (and hence has the same decisions as the optimal policy for $\bv{\Theta}_1$). 
Then: 
\[ \bv{J}_{\bv{\Theta}_2} \leq \bv{J}_{21}^{mis} \leq \bv{J}_{\bv{\Theta}_1} + \bv{1}\frac{\theta}{\phi} \]
where $\bv{1}$ is a vector of all $1$ values with the same dimension as $\bv{J}_{\bv{\Theta}_1}$.

\begin{proof} 
By definition, we have $\bv{J}_{\bv{\Theta}_2} \leq \bv{J}_{21}^{mis}$ (as $\bv{J}_{\bv{\Theta}_2}$ is the 
minimum sum cost over any policy when penalties are incurred according to $\bv{\Theta}_2$ queue backlog).   
Consider  any entry $z$, and suppose
we start in initial state $z(0)=z$.\footnote{Note that while all frames start with $z=0$, and hence have optimal
cost $J^*[0]$, $\bv{J}^*$ is defined with entries indexed by general initial states $z$.}   
  Then:
\begin{eqnarray*}
J_{\bv{\Theta}_2}[z] &\leq& J_{21}^{mis}[z] \\
&=&  \expect{\sum_{\tau=0}^{T-1} c_{\bv{\Theta}_2}(\alpha_1(\tau), \Omega(\tau), z_1(\tau))} \label{eq:diff2} \\
&=&  J_{\bv{\Theta}_1}[z]  +  \expect{\sum_{\tau=0}^{T-1} c_{\bv{\Theta}_2}(\alpha_1(\tau), \Omega(\tau), z_1(\tau))} \nonumber \\
&&-\expect{\sum_{\tau=0}^{T-1} c_{\bv{\Theta}_1}(\alpha_1(\tau), \Omega(\tau), z_1(\tau))} \label{eq:diff3} \\
&\leq& J_{\bv{\Theta}_1}[z]  + \frac{\theta}{\phi} \label{eq:diff4} 
\end{eqnarray*}
where the final inequality is due to the fact that the mean renewal time is $1/\phi$, and from the 
fact that the $\theta$ value in \eqref{eq:beta-theta} bounds the difference in the $c_{\bv{\Theta}_1}(\cdot)$ and 
$c_{\bv{\Theta}_2}(\cdot)$ components.
This proves part (b).  

To prove part (a), note that part (b) implies: 
\[ \bv{J}_{\bv{\Theta}_2} \leq \bv{J}_{\bv{\Theta}_1} + \bv{1} \frac{\theta}{\phi} \]
However, switching the roles of $\bv{\Theta}_1$ and $\bv{\Theta}_2$, we can similarly derive
$\bv{J}_{\bv{\Theta}_1} \leq \bv{J}_{\bv{\Theta_2}} + \bv{1}\theta/\phi$. This proves part (a). 
\end{proof}

\bibliographystyle{unsrt}
\bibliography{../../../../../../latex-mit/bibliography/refs}
\end{document}